# SYMMETRIC GIBBS MEASURES

## KARL PETERSEN AND KLAUS SCHMIDT

ABSTRACT. We prove that certain Gibbs measures on subshifts of finite type are nonsingular and ergodic for certain countable equivalence relations, including the orbit relation of the adic transformation (the same as equality after a permutation of finitely many coordinates). The relations we consider are defined by cocycles taking values in groups, including some nonabelian ones. This generalizes (half of) the identification of the invariant ergodic probability measures for the Pascal adic transformation as exactly the Bernoulli measures—a version of de Finetti's Theorem. Generalizing the other half, we characterize the measures on subshifts of finite type that are invariant under both the adic and the shift as the Gibbs measures whose potential functions depend on only a single coordinate. There are connections with and implications for exchangeability, ratio limit theorems for transient Markov chains, interval splitting procedures, 'canonical' Gibbs states, and the triviality of remote sigma-fields finer than the usual tail field.

## 1. INTRODUCTION

It has long been known that many Gibbs measures on (topologically mixing) subshifts of finite type (SFT's) are not only ergodic but are, in fact, $K$ (they satisfy the Kolmogorov 0-1-Law in that they have trivial (one-sided) tail fields). In fact the two-sided tail field is also trivial; another way to state this is to say that they are ergodic for the *homoclinic* or *Gibbs* relation (the countable equivalence relation $R_A$ in which two sequences are equivalent if and only if they disagree in only finitely many coordinates, i.e. if and only if they are in the same orbit under the action of the group $\Gamma$ of finite coordinate changes). We prove that these Gibbs measures, which include the Markov measures fully supported on the subshift, are also nonsingular and ergodic for many *cocycle-generated* subrelations of the homoclinic relation, including the

---

1991 *Mathematics Subject Classification.* Primary: 28D05, 60G09. Secondary: 58F03, 60J05, 60K35, 82B05.

First author supported in part by NSF Grant DMS-9203489.

*symmetric* relation $S_A$ in which two sequences are equivalent if and only if one can be obtained from the other by a *permutation* of finitely many coordinates. On a one-sided SFT this orbit relation for the action of the group $\Pi$ of permutations of finitely many coordinates is the orbit relation of the *adic transformation*, and ergodicity of this equivalence relation is the same as ergodicity of the adic transformation.

Our interest in these matters arose from the study of the dynamics of adic transformations, which were defined by A.M. Vershik (see [56]) as a family of models in which the cutting and stacking constructions of ergodic theory are organized in a way that makes them conveniently accessible to combinatorial analysis, and which has certain universality properties (such as containing a uniquely ergodic version of every ergodic, measure-preserving transformation on a Lebesgue space [55]). 'Transversal flows' had been investigated by Sinai [49], Kubo [29] and Kowada [26, 27], and also by S. Ito [24], who proved ergodicity of the adic transformation on an SFT for its measure of maximal entropy, a particular case of our Theorem 3.3. Adic transformations for actions of amenable groups were constructed by Lodkin and Vershik [35]. For the Pascal adic transformation (first defined by Vershik [54]), i.e. the adic transformation on the full shift, the ergodicity of Gibbs measures (including Bernoulli and Markov measures) is closely related to classical investigations of exchangeability and the triviality of remote sigma-fields. Given a (usually finite-valued) stochastic process $X_0, X_1, \ldots$ with (usually shift-invariant) distribution $\mu$, we have the sigma-algebra $\mathcal{I}$ of *shift-invariant sets*, the *tail field* $\mathcal{F}_\infty^+ = \cap_{n \geq 0} \mathcal{B}(X_n, X_{n+1}, \ldots)$, and the *exchangeable sigma-algebra* $\mathcal{E}$ consisting of all sets invariant under the group $\Pi$ of permutations of $\mathbb{N}$ that move only finitely many coordinates. (As sub-sigma-algebras of the product sigma-algebra of countably many copies of a fixed measure algebra, $\mathcal{I} \subset \mathcal{F}_\infty^+ \subset \mathcal{E}$, and the inclusions can be proper. In the two-sided case there are also $\mathcal{F}_\infty^-$, the sigma-algebra generated by $\mathcal{F}_\infty^+$ and $\mathcal{F}_\infty^-$, and the two-sided tail $\mathcal{F}_\infty = \cap_{n \geq 1} \mathcal{B}(X_i : |i| \geq n)$.) Already for the case of the full shift, invariance and ergodicity under the symmetric subrelation $S_n^+$ of the Gibbs relation $R_n^+$ of $\Sigma_n^+$ leads to interesting results. A. M. Vershik [private communication] observed that the nonatomic, ergodic, invariant probability measures for the Pascal adic transformation on the 2-shift are in one-to-one correspondence with the Bernoulli measures on $\Sigma_2^+$. By using dyadic expansions to regard this Pascal adic transformation as an infinite interval-exchange map on the unit interval $[0, 1]$, one may recognize it as the map studied by Arnold [3] and proved by Hajian, Ito, and Kakutani [18] to be ergodic (with respect to Lebesgue measure) and later used by Kakutani [25] to prove the uniform distribution of his interval splitting procedure. (The connection was also noted in [50]. See [42] for more about interval splitting.) The Hewitt-Savage 0-1-Law [20] says that for i.i.d. random variables $\mathcal{E}$ is trivial, i.e. consists only of sets of measure 0 or 1 (Bernoulli measures are $S_n^+$-ergodic). This was extended to Markov measures by Blackwell and Freedman [5]. There are analogous statements for the two-sided case $\ldots X_{-1}, X_0, X_1, X_2, \ldots$, where there are past, future, and two-sided tail fields to deal with. It took some time to sort out the relationships among the various tail fields [4, 17, 22, 23, 33, 37, 38, 40, 39], to understand the representation of measures as mixtures in terms of the theory of Choquet simplices and ergodic decompositions [9, 8, 10, 43, 45, 57], and to investigate exchangeability in more general contexts [9, 8, 14, 21, 31, 43, 52]. For a survey of exchangeability, see [2]. Diaconis and Freedman [8] gave a necessary and sufficient condition for a measure to be a mixture of Markov measures, in terms of 'partial exchangeability'—invariance under the subgroup of Π that preserves transition counts (as well as symbol counts). They also gave a general theory of 'sufficient statistics', describing how to present the (in some sense) most general symmetric measure as a mixture of extremal ones [9]. Several workers in statistical mechanics considered 'canonical' or 'microcanonical' Gibbs states, in the construction of which symbol counts, or the values of some other 'energy function', are fixed in finite regions. The point was to find the extremal measures in at least some classes of examples and thereby obtain a mixture-representation theorem for the most general such measures. Georgii [14, 15] found that certain Markov measures qualify. Similar results were obtained by Lauritzen [31], Höglund [21], and Thompson [52].

We prove (Theorem 3.3) that many Gibbs measures, including all mixing Markov measures, are *ergodic* under certain subrelations of $R_A$ which include the equivalence relation $S_A$ that corresponds to permutation of finitely many coordinates (or, equivalently, eventual equality of accumulated symbol counts). As a corollary, (letting $k = 2$ for simplicity) taking $\phi(x) = \psi(x) = cx_0$ in Theorem 3.3, we find measures $\mu_\phi$ that are ergodic and invariant under the adic on the SFT $\Sigma_A$, i.e. symmetric or exchangeable, when sequences are constrained to lie in the SFT $\Sigma_A$ (see Examples 5.1, 5.3, and 5.4). Similar results hold when $\psi$ is a function of only finitely many coordinates of $x$. Using a well-known formula (see [41, p. 22]), one can recover explicit expressions for these Gibbs states (cf. [14, 15, 21, 31, 52]). As a consequence of the ergodicity of Markov measures for $S_A^+$ we obtain some pointwise ratio limit theorems for adic transformations whose direct combinatorial proofs appear to be quite difficult (Section 7.4). Further, we identify all the probability measures on SFT's that are invariant simultaneously for the shift and the symmetric relation $S_A$ (Theorem 6.2). More generally, we prove ergodicity of certain Gibbs measures under subrelations $S_A^\psi$ of $R_A$ defined by cocycles generated by functions $\psi : \Sigma_A \longmapsto \mathcal{G}$, where $\mathcal{G}$ is a 'nearly abelian' countable discrete group—see 3.13 and Theorem 3.3; the symmetric relation $S_A$ is of this kind for a particular $\psi$. The dynamical and abstract approach provides a common framework for the diverse questions and results scattered through the literature and also leads to some interesting problems, mentioned in the final section of the paper, for example whether these measures are weakly mixing and whether there are some implications for various kinds of interval splitting or for the statistical mechanics of materials.

The authors thank C. Ji, U. Krengel, and A. M. Vershik for their suggestions that contributed significantly to the progress of this work and the referee for several improvements in the presentation.

## 2. Background

In this section we review the elements that we will need from the general theory of Borel equivalence relations, define adic transformations, and recall what is known about ergodic invariant measures for the Pascal adic transformation on the 2-shift.

2.1. **Borel equivalence relations.** First we establish the basic terminology and notation concerning Borel equivalence relations that will be needed later (cf. [11, 12, 48]).

Let $(X, \mathcal{B})$ be a standard Borel space, and let $R \subset X \times X$ be a discrete Borel equivalence relation, i.e. a Borel subset which is an equivalence relation, and which satisfies in addition that the *equivalence class*

$$R(x) = \{x' \in X : (x, x') \in R\}$$

of every point $x \in X$ is countable. Under this hypothesis the *saturation*

$$R(B) = \bigcup_{x \in B} R(x)$$

of every Borel set $B \in \mathcal{B}$ is again a Borel set. The *full group* $[R]$ of $R$ is the group of all Borel automorphisms $W$ of $X$ with $Wx \in R(x)$ for every $x \in X$. There exists a countable subgroup $G \subset [R]$ such that

$$Gx = \{gx : g \in G\} = R(x)$$

for every $x \in X$. A sigma-finite measure $\mu$ on $\mathcal{B}$ is *quasi-invariant* under $R$ if $\mu(R(B)) = 0$ for every $B \in \mathcal{B}$ with $\mu(B) = 0$; it is *ergodic* if, in addition, either $\mu(R(B)) = 0$ or $\mu(X \smallsetminus R(B)) = 0$ for every $B \in \mathcal{B}$.

Let $\mu$ be a probability measure on $\mathcal{B}$ which is quasi-invariant under $R$. Then $\mu$ is also quasi-invariant under every $W \in [R]$. In particular, if $G \subset [R]$ is a countable subgroup with $Gx = R(x)$ for every $x \in X$, then one can patch together the Radon-Nikodym derivatives $d\mu g/d\mu$, $g \in G$, and define a Borel map $\rho_\mu \colon R \longmapsto (0, \infty) \subset \mathbb{R}$ with the following properties:

(1) $\rho_\mu(Wx, x) = (d\mu W/d\mu)(x)$ $\mu$-a.e., for every $W \in [R]$,
(2) $\rho_\mu(x, x')\rho_\mu(x', x'') = \rho_\mu(x, x'')$ whenever $(x, x'), (x, x'') \in R$.

The map $\rho_\mu$ is the *Radon-Nikodym derivative* of $\mu$ under $R$, and $\mu$ is $(R\text{-})invariant$ if $\rho_\mu(x, x') = 1$ for all $(x, x') \in R \smallsetminus ((X \times N) \cup (N \times X))$, where $N \in \mathcal{B}$ is a $\mu$-null set.

2.2. **The main examples.** In order to describe a class of adic transformations of particular interest to us we assume that $V = (V_k, \, k \geq 0)$ is a sequence of finite, nonempty sets and put $X_V = \prod_{k \geq 0} V_k$, with the product of the discrete topologies. Two elements $x = (x_k)$, $x' = (x'_k)$ in $X_V$ are *comparable* (in symbols: $x \sim x'$) if $x_k = x'_k$ for all sufficiently large $k \geq 0$. The set

$$R_V = \{(x, x') \in X_V \times X_V : x \sim x'\} \tag{2.1}$$

is a Borel set and an equivalence relation, and each equivalence class

$$R_V(x) = \{x' \in X_V : x' \sim x\} \tag{2.2}$$

of $R_V$ is countable. In other words, $R_V$ is a discrete, standard equivalence relation in the sense of [11]; it is called the *homoclinic equivalence relation* or *Gibbs relation* of $X_V$. More generally, if $Y \subset X_V$ is a nonempty closed set, we denote by

$$R_Y = R_V \cap (Y \times Y) \tag{2.3}$$

the *Gibbs relation* of $Y$ and observe that $R_Y(y) = R_V(y) \cap Y$ and $R_Y(B) = R_V(B) \cap Y$ for every $y \in Y$ and $B \subset Y$. The following conditions on $Y$ are easily seen to be equivalent:

(T1) the equivalence relation $R_Y$ on $Y$ is topologically transitive, i.e. there exists a $y \in Y$ whose equivalence class $R_Y(y)$ is dense in $Y$,
(T2) for every pair $\mathcal{O}_1, \mathcal{O}_2$ of nonempty open subsets of $Y$ there exists a point $y \in Y$ with $R_Y(y) \cap \mathcal{O}_i \neq \varnothing$ for $i = 1, 2$.

If $Y$ fulfills (either of) these conditions we shall simply say that $Y$ *satisfies condition* (T).

Now assume that each of the sets $V_k$ is totally ordered with an order $<_k$. The sequence $<$ of orders $(<_k, \, k \geq 0)$ on the sets $V = (V_k, \, k \geq 0)$ induces a total order $\prec$ on each equivalence class $R_V(x)$, $x \in X_V$, and thus a partial order on the space $X_V$: if $x', x'' \in R_V(x)$ then $x' \prec x''$ whenever $x'_n <_n x''_n$ with $n = \max\{k \geq 0 : x'_k \neq x''_k\}$. In this partial order $X_V$ has a unique maximal element $x^+$ and a unique minimal element $x^-$; furthermore, if $x \in X_V$, $x \neq x^+$, then there exists a unique smallest element $succ(x) \in \{x' \in R_V(x) : x \prec x'\}$, called the *successor of* $x$. Put

$$C_V = R_V(x^-) \cup R_V(x^+).$$

The restriction of the map $x \mapsto succ(x)$ to $X_V \smallsetminus C_V$ is a homeomorphism of $X_V \smallsetminus C_V$; in order to extend this homeomorphism to a Borel automorphism of the entire space $X_V$ we set

$$T_V^{\prec}(x) = \begin{cases} succ(x) & \text{if } x \in X_V \smallsetminus C_V, \\ x & \text{if } x \in C_V, \end{cases}$$

and call $T_V = T_V^{\prec}$ the *adic transformation of* $X_V$ *determined by the orders* $\prec$ *of the sets* $V$.

In contrast to the space $X_V$, a nonempty, closed subset $Y \subset X_V$ may have many minimal and maximal elements with respect to the partial order $\prec$; however, if $y \in Y$ is not maximal, then there still exists a unique smallest element $succ(y) \in \{y' \in R_V(y) \cap Y : y \prec y'\}$. If $R_Y = R_V \cap (Y \times Y)$, and if $Y^-$ and $Y^+$ are the sets of minimal and maximal elements of $Y$ and $C_Y = R_Y(Y^-) \cup R_Y(Y^-)$, then the map $y \mapsto succ(y)$ again induces a Borel automorphism $T_Y : Y \longmapsto Y$, given by

$$T_Y y = \begin{cases} succ(y) & \text{if } y \in Y \smallsetminus C_Y, \\ y & \text{if } y \in C_Y, \end{cases}$$

which is called the *adic transformation* of $Y$. Note that

$$\{T_Y^n(y) : n \in \mathbb{Z}\} = \{y' \in Y : y' \sim y\} \tag{2.4}$$

for every $y \in Y \smallsetminus C_Y$.

We could have tried to define the adic transformation $T_Y$ a little more elegantly on the set $C_Y$ of exceptional points (cf. [16]). However, the definition of $T_Y$ on $C_Y$ will not really matter, since we shall only consider nonatomic measures and adic transformations $T_Y$ on closed, nonempty subsets $Y \subset X_V$ which satisfy the following condition:

(M)  the sets $Y^-$ and $Y^+$ of minimal and maximal elements of $Y$ in the partial order $\prec$ are both countable.

*Example 2.1.* (*The Pascal adic transformation (Vershik* [54]*)*) Let $n \geq 2$, and let

$$V_k^{(n)} = \{(j^{(0)}, \dots, j^{(n-1)}) \in \mathbb{N}^n : j^{(0)} + \cdots + j^{(n-1)} = k\}$$

for every $k \geq 0$. We put $V^{(n)} = (V_k^{(n)}, \ k \geq 0)$, furnish each $V_k^{(n)}$ with the reverse lexicographic order, and set

$$Y^{(n)} = \{y = (y_k) \in X_{V^{(n)}} : y_k = (y_k^{(0)}, \ldots, y_k^{(n-1)}) \in V_k^{(n)} \ \text{and}$$

$$y_k^{(i)} \leq y_{k+1}^{(i)} \ \text{for every} \ k \geq 0 \ \text{and} \ i = 0, \ldots, n-1\}. \tag{2.5}$$

The equivalence relation $R_{Y^{(n)}}$ and the adic transformation $T_n = T_{Y^{(n)}}$ are called the $n$-dimensional Pascal Gibbs relation and $n$-dimensional Pascal adic transformation on the space $Y^{(n)}$. It is easy to see that $Y^{(n)}$ satisfies the conditions (T) and (M).

We picture $Y^{(n)}$ as a graded graph. At level $k$ there are as vertices the elements of $V_k^{(n)}$. There is a connection from $y_k \in V_k^{(n)}$ to $y_{k+1} \in V_{k+1}^{(n)}$ if and only if $y_k^{(i)} \leq y_{k+1}^{(i)}$ for all $i = 0, \ldots, n-1$. In fact then $y_k^{(i)} = y_{k+1}^{(i)}$ for all $i$ except one, $i = i_0$, for which $y_k^{i_0} = y_{k+1}^{i_0} - 1$. We think of the edge from $y_k$ to $y_{k+1}$ as labeled with the symbol $i_0 \in \{0, \ldots, n-1\}$; denote this $i_0$ by $\phi_k(y)$. Then $y_k^{(j)}$ is equal to the number of appearances of symbol $j$ on the edges of the path $y$ from level $0$ to level $k$.

This allows us to define a continuous bijection

$$\Phi \colon Y^{(n)} \longmapsto \Sigma_n^+ = (\mathbb{Z}/n\mathbb{Z})^{\mathbb{N}}$$

by setting, for every $y = (y_k) \in Y^{(n)}$,

$$\Phi(y)_k = \phi_k(y).$$

If we put $W^{(n)} = (W_k^{(n)}, \ k \geq 0)$ with $W_k^{(n)} = \mathbb{Z}/n\mathbb{Z}$ for every $k \geq 0$, and if we furnish each $W_k^{(n)}$ with its natural order, then $\Sigma_n^+ = X_{W^{(n)}}$, and we define the partial order $\prec'$ on $\Sigma_n^+ = X_{W^{(n)}}$ as above. It is clear that the map $\Phi \colon Y^{(n)} \longmapsto \Sigma_n^+$ is order-preserving, i.e. that $\Phi(y) \prec' \Phi(y')$ if and only if $y \prec y'$. However, if $R_n^+ = R_{W^{(n)}}$ is the Gibbs relation on $\Sigma_n^+$, then

$$(\Phi \times \Phi)(R_{Y^{(n)}}) = S_n^+ \subsetneq R_n^+. \tag{2.6}$$

Indeed, $S_n^+$ is the subrelation of $R_n^+$ in which two points $(z, z') \in R_n^+$ are equivalent if and only if the coordinates of $z'$ differ from those of $z$ by a finite permutation.

There exists a unique Borel automorphism $T_{\Sigma_n^+}^*$ of $\Sigma_n^+$ satisfying

$$T_{\Sigma_n^+}^* \cdot \Phi(y) = \Phi \cdot T_{Y^{(n)}}$$

for every $y \in Y^{(n)}$. The map $T_{\Sigma_n^+}^*$ satisfies that

$$\{(T_{\Sigma_n^+}^*)^k(x) : k \in \mathbb{Z}\} = S_n^+(x)$$

for all but countably many $x \in \Sigma_n^+$. It is important to note the distinction between the map $T_{\Sigma_n^+}^*$ just defined and the adic transformation $T_{\Sigma_n^+}$ arising from the partial order $\prec'$ on $\Sigma_n^+ = X_{W^{(n)}}$, the familiar adding machine or odometer, most of whose orbits are equal to the equivalence classes of the Gibbs relation $R_n^+$; in the case $n = 2$, $T_{\Sigma_2^+}^*(0^p1^q10\ldots) = 1^q0^p01\ldots$ for each $p, q \geq 0$. We shall abuse terminology to some extent and call the transformation $T_{\Sigma_n^+}^*$ the *adic transformation on the full shift* $\Sigma_n^+$.

*Example 2.2.* (*The adic transformation of an SFT*) Define $V^{(n)} = (V_k^{(n)}, k \geq 0)$ and $X_{V^{(n)}}$ for $n \geq 1$ as in Example 2.1, and let $A = (A(i,j), 0 \leq i, j \leq n-1)$ be an irreducible, aperiodic transition matrix with entries in $\{0, 1\}$. We denote by

$$\Sigma_A^+ = \{x = (x_k) \in \Sigma_n^+ : A(x_k, x_{k+1}) = 1 \text{ for every } k \in \mathbb{N}\}$$

the one-sided SFT defined by $A$ and write

$$R_A^+ = R_n^+ \cap (\Sigma_A^+ \times \Sigma_A^+) \tag{2.7}$$

for the Gibbs relation of $\Sigma_A^+$. Put

$$Y_A = \Phi^{-1}(\Sigma_A^+) \subset X_{V^{(n)}}. \tag{2.8}$$

Then $Y_A$ satisfies the conditions (T) and (M) above (since any maximal path in the associated graph, when traced backwards, must always follow the allowed edge with the largest possible label, and hence the labelling has to be eventually periodic).

Let $S_A^+ \subset R_A^+$ be the equivalence relation in which two points $(z, z') \in R_A^+$ are equivalent if there is $K$ such that $z_k = z_{k'}$ for $k \geq K$ while for $k < K$ the coordinates $z_k$ of $z$ are permutations of those of $z'$. Then $(\Phi \times \Phi)(R_{Y_A}) = S_A^+$. The Borel automorphism $T_{\Sigma_A^+}^* : \Sigma_A^+ \smallsetminus C_{\Sigma_A^+} \longmapsto \Sigma_A^+ \smallsetminus C_{\Sigma_A^+}$ satisfying $T_{\Sigma_A^+}^* \cdot \Phi(y) = \Phi \cdot T_{Y_A}(y)$ for every $y \in Y_A \smallsetminus C_A$ can again be determined explicitly. For example, if $n = 2$ and $A = \left(\begin{smallmatrix} 1 & 1 \\ 1 & 0 \end{smallmatrix}\right)$, then $\Sigma_A^+ \subset \Sigma_2^+$ is the *golden mean SFT*, and $T_{\Sigma_A^+}^*(0^p(10)^q100\ldots) = (10)^q0^p010\ldots$ for every $p \geq 0$, $q \geq 1$.

In analogy with Example 2.1 we call the transformation $T_{\Sigma_A^+}^*$ the *adic transformation of the SFT* $\Sigma_A^+$, but emphasize that $T_{\Sigma_A^+}^*$ *differs from the transformation* $T_{\Sigma_A^+}$, the corresponding *stationary adic* transformation (in the case of the full shift, the odometer), whose orbits generate the equivalence classes of $R_A^+$ for all but countably many points in $\Sigma_A^+$.

*Example 2.3.* (*The Gibbs relation and adic transformation of a graph*) Let $V = (V_k, k \geq 0)$ be a sequence of finite, nonempty sets, and let, for every $k \geq 0$, $E_k \subset V_k \times V_{k+1}$ be a subset whose projection on each of the two sets $V_k$ and $V_{k+1}$ is surjective. We set $E = (E_k, \ k \geq 0)$ and put

$$Y_{V,E} = \{x = (x_k) \in X_V : [x_k, x_{k+1}] \in E_k \text{ for every } \ k \geq 0\}. \tag{2.9}$$

If $\bar{V}_k = \{(k,v) : v \in V_k\}$, $\bar{E}_k = \{(k,v) \to (k+1,v') : [v,v'] \in E_k\}$, $\bar{V} = \bigcup_{k \geq 0} \bar{V}_k$, $\bar{E} = \bigcup_{k \geq 0} \bar{E}_k$, and if $\Gamma = \Gamma_{V,E}$ is the directed graph with vertices $\bar{V}$ and edges $\bar{E}$, then there is a bijective correspondence between the elements of $Y_{V,E}$ and the infinite paths in $\Gamma$ which start at an element of $\bar{V}_0$. The space $Y_{V,E}$ satisfies the conditions (T) if and only if there exist, for every pair of vertices $\bar{v}, \bar{v}' \in \bar{V}$, a vertex $\bar{w} \in \bar{V}$ and paths $\bar{v} \to \bar{v}_1 \to \cdots \to \bar{v}_m \to \bar{w}$, $\bar{v}' \to \bar{v}'_1 \to \cdots \to \bar{v}'_{m'} \to \bar{w}$. The relation $R_{V,E}$ is called the *Gibbs relation* of the directed graph $\Gamma$; if the sets $V_k$, $k \geq 0$, are totally ordered, and if $Y = Y_{V,E}$ satisfies (M) then the adic transformation $T_Y$ is the *adic transformation of the graph* $\Gamma$.

In general, if $Y \subset X_V$ is a nonempty, closed subset, we denote by $\pi_S \colon Y \longmapsto \prod_{k \in S} V_k$ the restriction to $Y$ of the coordinate projection $\prod_{k \geq 0} V_k \longmapsto \prod_{k \in S} V_k$. Then $Y$ is of the form $Y = Y_{V,E}$ described above if and only if

$$\pi_{\{0,\dots,n\}}(R_Y^{(n)}(y)) = \pi_{\{0,\dots,n\}}(R_Y^{(n)}(y')) \tag{2.10}$$

whenever $y = (y_k), y' = (y'_k) \in Y$ and $y_n = y'_n$, where $R_Y^{(n)}$ is defined in (2.12). The Gibbs relation in Example 1.1 satisfies (2.10) and is thus the Gibbs relation of a graph, but the relation $R_A$ in Example 1.2 need not satisfy (2.10).

2.3. **Adic-invariant sets and measures.** Throughout this section we assume that $V = (V_k, \ k \geq 0)$ is a sequence of finite, nonempty, totally ordered sets, put $X_V = \prod_{k \geq 0} V_k$, define the equivalence relation $R_V$ as in Section 2.1, and assume that $Y \subset X_V$ is a closed subset and $R_Y = R \cap (Y \times Y)$.

Let $\mathcal{B}_Y$ be the Borel field of $Y$, and let $\mu$ be a nonatomic probability measure on $\mathcal{B}_Y$. If $Y$ satisfies the condition (M), and if $T_Y$ is the adic transformation on $Y$, then $\mu$ is quasi-invariant (or ergodic) under $R_Y$ if and only if it is quasi-invariant (or ergodic) under $T_Y$, and the sigma-algebras

$$\begin{aligned}
\mathcal{B}_Y^{R_Y} &= \{R_Y(B) : B \in \mathcal{B}_Y\}, \\
\mathcal{B}_Y^{T_Y} &= \{B \in \mathcal{B}_Y : T_Y B = B\}
\end{aligned} \tag{2.11}$$

of $R_Y$-saturated and $T_Y$-invariant Borel sets in $Y$ coincide (mod $\mu$).

We denote by $\mathcal{T}_n \subset \mathcal{B}_Y$ the sigma-field generated by all cylinder sets of the form $[v_n, \ldots, v_{n+k}]_n = \{y = (y_i) \in Y : y_n = v_n, \ldots, y_{n+k} = v_{n+k}\}$ and define the *tail sigma-field* of $Y$ as $\mathcal{T}_Y = \bigcap_{n \geq 0} \mathcal{T}_n$.

**Lemma 2.4.** *For every closed, nonempty subset* $Y \subset X_V$, $\mathcal{T}_Y = \mathcal{B}_Y^{R_Y}$. *In particular, if* $\mu$ *is a nonatomic probability measure on* $Y$ *which is quasi-invariant under* $R_Y$ *(or, if* $Y$ *satisfies the hypothesis* (M) *under the adic transformation* $T_Y$*), then*

$$\mathcal{T}_Y = \mathcal{B}_Y^{R_Y} = \mathcal{B}_Y^{T_Y} \pmod{\mu}.$$

*Proof.* For every $m \geq 0$ we set

$$R_Y^{(m)} = \{(y, y') \in Y \times Y : y_k = y'_k \text{ for all } k \geq m\} \tag{2.12}$$

and note that $R_Y = \bigcup_{n \geq 0} R_Y^{(n)}$ and hence $\mathcal{B}_Y^{R_Y} = \bigcap_{n \geq 0} \mathcal{B}^{R_Y^{(n)}}$. As $\mathcal{T}_n = \mathcal{B}_Y^{R_Y^{(n)}}$ for every $n \geq 0$, we conclude that $\mathcal{T}_Y = \mathcal{B}_Y^{R_Y}$. $\qquad\square$

The *m-weight* $w_m(y)$ of an element $y \in Y \subset X_V$, $m \geq 0$, is defined as

$$w_m(y) = |R_Y^{(m)}(y)| \tag{2.13}$$

where $|F|$ denotes the cardinality of a set $F$. Similarly, if $C \subset Y$ is a Borel set and $m \geq 0$, then

$$w_m(C, y) = |R_Y^{(m)}(y) \cap C|. \tag{2.14}$$

(This corresponds to the *dimension* of a vertex in a Bratteli diagram, the number of paths into a vertex in a graph as in Example 2.3, or the height of a column in an ergodic-theoretic cutting and stacking construction.) For example, if $Y^{(n)} \subset X_{V^{(n)}}$ is defined as in Example 2.1, then the $m$-weight $w_m(y)$ of an element $y = (y_k) \in Y^{(n)}$ is given as follows: if $y_m = (y_m^{(0)}, \ldots, y_m^{(n-1)}) \in V_m^{(n)}$, then

$$w_m(y) = \frac{m!}{y_m^{(0)}! \cdots y_m^{(n-1)}!}.$$

**Theorem 2.5** (Vershik [53])**.** *Let* $Y \subset X_V$ *be a closed, nonempty subset, and let* $\mu$ *be a nonatomic probability measure on* $\mathcal{B}_Y$ *which is invariant and ergodic under* $R_Y$*. For every Borel set* $B \subset Y$,

$$\mu(B) = \lim_{n \to \infty} \frac{w_n(B, y)}{w_n(y)} \qquad \text{for } \mu\text{-a.e. } y \in Y,$$

*where the weights* $w_n(\cdot)$ *and* $w_n(B, \cdot)$ *are defined in* (2.13) *and* (2.14)*.*

*Proof.* If $f \in L^1(Y, \mathcal{B}_Y, \mu)$, then

$$E_\mu(f|\mathcal{T}_n)(y) = E_\mu(f|\mathcal{B}_Y^{R_Y^{(n)}})(y) = \frac{1}{|R_Y^{(n)}|} \sum_{z \in R_Y^{(n)}} f(z)$$

for $\mu$-a.e. $y \in Y$. The reverse martingale theorem and Lemma 2.4 imply that

$$\lim_{n \to \infty} \frac{1}{|R_Y^{(n)}(y)|} \sum_{z \in R_Y^{(n)}(y)} f(z) = E_\mu(f|\mathcal{T}_Y)(y) = E_\mu(f|\mathcal{B}_Y^{R_Y})(y)$$

$\mu$-a.e., and by setting $f$ equal to the indicator function $1_B$ of $B$ we have proved the theorem. $\qed$

*Remark 2.6.* If $\mu$ is invariant but not necessarily ergodic, then

$$\frac{w_n(B, y)}{w_n(y)} \to E_\mu(1_B|\mathcal{T})(y) \qquad \mu\text{-a.e.}$$

for every Borel set $B \subset Y$.

**Theorem 2.7** (Hajian-Ito-Kakutani [18]). *Let $Y^{(2)}$ be as in (2.5), and denote by $\Phi \colon Y^{(2)} \longmapsto \Sigma_2^+ = \{0, 1\}^{\mathbb{N}}$ the map described in* Example 2.1, *so that $S_2^+ = (\Phi \times \Phi)(R_{Y^{(2)}})$ is the equivalence relation on $\Sigma_2^+$ in which two points are equivalent if and only if their coordinates differ by a finite permutation. For any $\alpha$ with $0 < \alpha < 1$ we set $\nu_\alpha(0) = \alpha$, $\nu_\alpha(1) = 1 - \alpha$, write $\mu_\alpha = \nu_\alpha^{\mathbb{N}}$ for the corresponding Bernoulli measure on $\Sigma_2^+$, and put $\bar\mu_\alpha = \mu_\alpha \cdot \Phi$. Then the following equivalent statements are true:*

(1) *$\bar\mu_\alpha$ is invariant and ergodic under the Pascal adic transformation $T_2$,*

(2) *$\bar\mu_\alpha$ is invariant and ergodic under $R_{Y_2}$,*

(3) *$\mu_\alpha$ is invariant and ergodic under the adic transformation $T_{\Sigma_2^+}^*$ of the full shift $\Sigma_2^+$ (or, equivalently, under the equivalence relation $S_2^+$).*

*Proof.* The equivalence of (1), (2), and (3) is obvious; we prove (2). We write $C(m, n)$ for the binomial coefficient $\binom{m}{n} = m!/(n!(m-n)!)$. For each cylinder set $C = [v_0, \dots, v_m] = \{y \in Y_2 : y_i = v_i$ for $i = 0, \dots, m\} \subset Y_2$ ending in a coordinate $v_m = (j_m^{(0)}, j_m^{(1)}) \in V_m^{(1)}$ with $j_m^{(0)} + j_m^{(1)} = m$, and for every $y = (y_k) \in Y_2$, $n > m$ and $y_n = (y_n^{(0)}, y_n^{(1)})$ with $j_n^{(i)} \geq j_m^{(i)}$, $i = 0, 1$,

$$\frac{w_n(C, y)}{w_n(y)} = \frac{C(n-m, y_n^{(0)} - j_m^{(0)})}{C(n, y_n^{(0)})} = \frac{(n-m)! \, y_n^{(0)}! \, y_n^{(1)}!}{(y_n^{(0)} - j_m^{(0)})! \, (y_n^{(1)} - j_m^{(1)})! \, n!}$$

$$= \frac{(y_n^{(0)} - j_m^{(0)} + 1) \cdots y_n^{(0)} \cdot (y_n^{(1)} - j_m^{(1)} + 1) \cdots y_n^{(1)}}{(n-m+1) \cdots n}$$

$$= \left(\frac{y_n^{(0)}}{n}\right)^{j_m^{(0)}} \left(\frac{y_n^{(1)}}{n}\right)^{j_m^{(1)}} \cdot \prod_{j=1}^{j_m^{(0)} - 1} \left(1 - \frac{j}{y_n^{(0)}}\right) \cdot \prod_{j=1}^{j_m^{(1)} - 1} \left(1 - \frac{j}{y_n^{(1)}}\right)$$

$$\cdot \prod_{j=1}^{m-1} \left(1 - \frac{j}{n}\right)^{-1}$$

$$\to \alpha^{j_m^{(0)}} (1 - \alpha)^{(j_m^{(1)})}$$

for $\bar{\mu}_\alpha$-a.e. $y \in Y_2$ as $n \to \infty$. In particular, $\bar{\mu}_\alpha(C)$ only depends on the last coordinate $v_m$ of $C$, which guarantees the invariance of $\bar{\mu}_\alpha$ under $R_{Y_2}$, and $\bar{\mu}_\alpha$ is ergodic because $\lim_{n \to \infty} w_n(C, \cdot)/w_n(\cdot)$ is constant $\bar{\mu}_\alpha$-a.e. (cf. Remark 2.6). $\qquad \square$

*Remark 2.8.* Since the measures $\bar{\mu}_\alpha$, $\alpha \in (0, 1)$, in Theorem 2.7 are obtained from the Bernoulli measures $\mu_\alpha$ on $\Sigma_2^+$ via the map $\Phi$, they are called the *Bernoulli measures* on $Y_2$. If we also consider the degenerate Bernoulli measures $\mu_0, \mu_1$ on $\Sigma_2^+$ defined as above, but with $\alpha \in \{0, 1\}$, then the measures $\bar{\mu}_0, \bar{\mu}_1$ are again invariant and trivially ergodic under $R_{Y_2}$: each of them is concentrated on an equivalence class of $R_{Y_2}$ consisting of a single point.

**Theorem 2.9** (de Finetti, Vershik)**.** *The only nonatomic probability measures*
*which are invariant and ergodic for the Pascal adic transformation $T_2$ are the measures $\bar{\mu}_\alpha$, $\alpha \in (0, 1)$, described in* Theorem 2.7.

*Proof.* Let $C = [v_0, \ldots, v_m] = \{y \in Y_2 : y_i = v_i \text{ for } i = 0, \ldots, m\}$ be a cylinder set with $v_m = (j_m^{(0)}, j_m^{(1)}) \in V_m^{(2)}$, let $v_{m+1} = (j_{m+1}^{(0)}, j_{m+1}^{(1)}) \in V_{m+1}^{(2)}$ with $j_{m+1}^{(i)} \geq j_m^{(i)}$ for $i = 0, 1$, and put $C' = [v_0, \ldots, v_{m+1}] \subset Y_2$ (we are using the same notation as in Example 2.1 and Theorem 2.7). If $\mu$ is a nonatomic, ergodic, invariant probability measure for $T_2$, then Theorems 2.5 and 2.7 imply that

$$\frac{\mu(C)}{\mu(C')} = \lim_{n \to \infty} \frac{w_n(C, y)/w_n(y)}{w_n(C', y)/w_n(y)}$$

$$= \lim_{n \to \infty} \frac{C(n-m, y_n^{(0)} - v_m^{(0)})}{C(n-m-1, y_n^{(0)} - v_{m+1}^{(0)})} = 1 - \lim_{n \to \infty} \frac{y_n^{(0)}}{n}.$$

By Theorem 2.5 this limit exists $\mu$-a.e. and is $\mu$-a.e. equal to a constant $\alpha \in (0,1)$, so we conclude that $\mu = \bar{\mu}_\alpha$. $\qquad\blacksquare$

*Remark 2.10.* In order to make the connection between Theorem 2.9 and de Finetti's theorem more explicit, let $n \geq 2$, and let $\bar{\mu}$ be a nonatomic, ergodic, invariant probability measure for the adic transformation $T_n$ on $Y^{(n)}$ (cf. Example 2.1).Then $\mu = \bar{\mu}\Phi^{-1}$ is invariant under the relation $S_n^+$ of finite coordinate changes, and de Finetti's theorem shows that $\mu$ is a Bernoulli measure on $\Sigma_n^+$. More generally, every nonatomic invariant probability measure for $T_n$ is a mixture of such Bernoulli measures.

If $\Sigma_A^+$ is an arbitrary, irreducible, aperiodic SFT, then the explicit determination of all $S_A^+$-invariant and ergodic probability measures is not so obvious. In Section 6 we will answer this question for shift-invariant measures. In some special cases, de Finetti's theorem already gives a complete answer.

**Theorem 2.11.** *Let $A = \left(\begin{smallmatrix} 1 & 1 \\ 1 & 0 \end{smallmatrix}\right)$. Then the only nonatomic probability measures on $\Sigma_A^+ \subset \Sigma_2^+$ which are ergodic and invariant under $S_A^+$ (or, equivalently, under the adic transformation $T_{\Sigma_A^+}^*$) are the Markov measures $\mu_\alpha'$, $\alpha \in (0,1)$, where $\mu_\alpha'$ has transition matrix $\left(\begin{smallmatrix} \alpha & 1-\alpha \\ 1 & 0 \end{smallmatrix}\right)$ and initial distribution $p(0) = \alpha$, $p(1) = 1 - \alpha$ (cf. (3.1)).*

*Proof.* We define a continuous, bijective map $\Psi \colon \Sigma_A^+ \longmapsto \{a, b\}^{\mathbb{N}} \cong \Sigma_2^+$ by replacing each occurrence of the string 10 in an element $x \in \Sigma_A^+$ with a $b$ and the remaining 0's with $a$'s. Although the map $\Psi$ is not shift-commuting, it satisfies that $(\Psi \times \Psi)(S_A^+) = S_2^+$; in particular, $\Psi$ sends the set of nonatomic $S_A^+$-invariant and ergodic probability measures on $\Sigma_A^+$ bijectively onto the set of nonatomic $R_2^+$-invariant, ergodic probability measures on $\Sigma_2^+$.

The latter measures are characterized by de Finetti's theorem as the Bernoulli measures on $\Sigma_2^+$, and the measures $\mu_\alpha'$, $\alpha \in (0,1)$, are just their images under $\Psi^{-1}$. $\qquad\blacksquare$

## 3. Gibbs measures and subrelations of Gibbs equivalence relations on two-sided shift spaces

As we saw in Theorems 2.5, 2.7 and 2.9, the ergodicity of Bernoulli measures under the equivalence relation $R_n^+$ on the full $n$-shift has dynamical implications for the $n$-dimensional Pascal Gibbs relation $R_{Y^{(n)}}$ and the $n$-dimensional Pascal adic transformation $T_{\Sigma_n^*}^*$ on $Y^{(n)}$. In this section we generalize our discussion of ergodicity of Bernoulli measures under the Pascal adic transformation by showing that certain Gibbs measures on SFT's are ergodic under a natural class of subrelations of the Gibbs relation of the subshift. However, since the discussion of ergodicity of these subrelations is in some sense more natural on two-sided shift spaces, we first discuss equivalence relations on two-sided SFT's before turning to one-sided SFT's in the next section.

Assume that $A = (A(i,j),\ 0 \le i, j \le n-1)$ is an irreducible, aperiodic $n \times n$ transition matrix with entries in $\{0,1\}$, denote by

$$\Sigma_A = \{x = (x_k) \in \Sigma_n = (\mathbb{Z}/n\mathbb{Z})^{\mathbb{Z}} :$$
$$A(x_k, x_{k+1}) = 1 \ \text{for every} \ k \in \mathbb{Z}\}$$

the two-sided SFT defined by $A$, and write $\sigma = \sigma_A$ for the shift

$$\sigma(x)_k = x_{k+1}$$

and

$$R_A = \{(x, x') \in \Sigma_A \times \Sigma_A : x_k \ne x'_k \ \text{for only finitely many} \ k \in \mathbb{Z}\}$$

for the Gibbs relation on $\Sigma_A$.

A one-step Markov measure $\mu_P$ on $\Sigma_A$ is determined by a stochastic matrix $P = (P(i,j),\ 0 \le i, j \le n-1)$ with $P(i,j) > 0$ if and only if $A(i,j) = 1$ (such a matrix is said to be *compatible with A*): if $\bar{p} = (\bar{p}(0), \dots, \bar{p}(n-1))$ is the unique probability vector with $\bar{p}P = \bar{p}$, then the measure of every cylinder set

$$C = [v_0, \dots, v_l]_m = \{x = (x_k) \in \Sigma_A : x_{m+i} = v_i \ \text{for} \ i = 0, \dots, l\}$$

is given by

$$\mu_P(C) = \bar{p}(v_0) P(v_0, v_1) \cdots P(v_{l-1}, v_l). \tag{3.1}$$

Any such Markov measure $\mu_P$ is easily seen to be quasi-invariant under the Gibbs relation $R_A$: its Radon-Nikodym derivative is given by

$$\rho_{\mu_P}(x, x') = \prod_{k=-\infty}^{\infty} \frac{P(x_k, x_{k+1})}{P(x'_k, x'_{k+1})} \tag{3.2}$$

for every $(x, x') \in R_A$ (note that the infinite product in (3.2) consists mostly of 1's).

More generally, let $\phi \colon \Sigma_A \longmapsto \mathbb{R}$ be a continuous function and put, for every $k \ge 0$,

$$\omega_0(\phi) = \max\{|\phi(x) - \phi(x')| : x, x' \in \Sigma_A\},$$
$$\omega_k(\phi) = \max\{|\phi(x) - \phi(x')| : x_l = x'_l \ \text{for} \ |l| < k\}, \ k \ge 1.$$

The function $\phi$ has *summable variation* if

$$\omega(\phi) = \sum_{k \geq 0} \omega_k(\phi) < \infty. \tag{3.3}$$

We denote by $M_1(\Sigma_A)$ the set of Borel probability measures on $\Sigma_A$, furnished with the weak* topology. A measure $\mu \in M_1(\Sigma_A)$ is a *Gibbs measure* of a map $\phi \colon \Sigma_A \longmapsto \mathbb{R}$ with summable variation if

$$\log \rho_\mu(x, x') = \sum_{k \in \mathbb{Z}} (\phi \cdot \sigma^k(x) - \phi \cdot \sigma^k(x')) \tag{3.4}$$

for every $(x, x') \in R_A \setminus (N \times N)$, where $N \subset \Sigma_A$ is a $\mu$-null set; the set of all such measures is denoted by $M_1^\phi(\Sigma_A)$. By definition, any Gibbs measure of $\phi$ is quasi-invariant under $R_A$. The Markov measure $\mu_P$ in (3.1) is a Gibbs measure of the map $\phi(x) = \log P(x_0, x_1)$. Conversely, if the function $\phi$ depends only on the coordinates $x_0, x_1$, then there is only one Gibbs measure for $\phi$, and it is Markov: $M_1^\phi(\Sigma_A) = \{\mu_P\}$ for some Markov matrix $P$ compatible with $A$ (see [41]). The following generalization of this fact is part of the lore of the theory of Gibbs measures; see, for example, [6, 32, 44, 58].

**Theorem 3.1.** *Let $A = (A(i, j), 0 \leq i, j \leq n-1)$ be an irreducible and aperiodic 0-1-matrix, $\Sigma_A \subset \Sigma_n$ the associated shift of finite type, and $R_A$ the Gibbs equivalence relation of $\Sigma_A$ described above. If $\phi \colon \Sigma_A \longmapsto \mathbb{R}$ is a function with summable variation, then there exists a unique Gibbs measure $\mu_\phi \in M_1^\phi(\Sigma_A)$. Moreover, $\mu_\phi$ is ergodic under $R_A$, and invariant and $K$ under the shift $\sigma = \sigma_A$.*

*Proof.* First we show that $M_1^\phi(\Sigma_A)$ is nonempty. For every $K \geq 0$ we set

$$R_A^{(K)} = \{(x, x') \in \Sigma_A \times \Sigma_A : x_k = x'_k \text{ for } |k| \geq K\}.$$

We claim that the set

$$M(\phi, K) = \{\mu \in M_1(\Sigma_A) : \mu \text{ satisfies (3.4) for every } (x, x') \in R_A^{(K)}\}$$

is nonempty for every $K \geq 0$. Indeed, fix $L > 1$ such that every entry of $A^L$ is positive, take $K \geq L$, and group the cylinder sets determined by central $(2K + 1)$-blocks into $n^2$ classes according to their first and last entries. We pick a member from each class and define $\mu$ on it arbitrarily, for example as a Markov measure consistent with the transition matrix $A$. Now we can use (3.4) to carry $\mu$ over to each of the other cylinder sets in the same class. Thus if $C$ is one of our chosen cylinder sets and $\gamma$ is the finite coordinate change that carries $C$ to $C'$ by changing its central $(2K + 1)$-block (of course the endpoints of the block do not change), then the summable variation property of $\phi$ shows that $\sum_{m\in\mathbb{Z}}[\phi\cdot\sigma^m(x)-\phi\cdot\sigma^m(\gamma x)]$ is a continuous function of $x\in C$, so for $E\subset C'$ we can define

$$\mu(E)=\int_{\gamma^{-1}E}\exp\sum_{m\in\mathbb{Z}}[\phi\cdot\sigma^m\gamma x)-\phi\cdot\sigma^m(x)]\,d\mu(x).$$

Since each $M(\phi,K)$ is also convex, compact, and nonincreasing in $K$, it follows that

$$M_1^\phi(\Sigma_A)=\bigcap_{K\geq 0}M(\phi,K)\neq\varnothing.$$

Next we show that $M_1^\phi(\Sigma_A)$ consists of a single measure which is ergodic under $R_A$; for this we require a little bit of notation and a lemma. Since the matrix $A$ is irreducible and aperiodic there exists an integer $L\geq 1$ such that $A^L(i,j)>0$ for every $(i,j)\in(\mathbb{Z}/n\mathbb{Z})^2$. A string $s=s_0\ldots s_r$ in $\{0,\ldots,n-1\}^{r+1}$ is *allowed* if $A(s_j,s_{j+1})=1$ for all $j=0,\ldots,r-1$. If $s'=s'_0\ldots s'_{r'}$ is a second allowed string, then $s$ and $s'$ can be *concatenated* if $s_r=s'_0$, and the string $ss'=s_0\ldots s_r s'_1\ldots s'_{r'}$ is the *concatenation* of $s$ and $s'$. Finally, if $s=s_0\ldots s_r$ is an allowed string and $m\in\mathbb{Z}$, then

$$[s]_m=\{x\in\Sigma_A:x_m=s_0,\ldots,x_{m+r}=s_r\}.$$

**Lemma 3.2.** *Suppose that* $\phi\colon\Sigma_A\longmapsto\mathbb{R}$ *has summable variation and* $\mu\in M_1^\phi(\Sigma_A)$. *Given* $M\geq 0$ *and disjoint cylinder sets*

$$C=[i_{-M}\ldots i_M]_{-M}=\{x\in\Sigma_A:x_l=i_l\ \text{for}\ |l|\leq M\}$$
$$D=[j_{-M}\ldots j_M]_{-M}=\{x\in\Sigma_A:x_l=j_l\ \text{for}\ |l|\leq M\},$$

*there is a homeomorphism* $V\in[R_A]$ *and a cylinder set* $C'\subset C$ *such that*

(1) $V^2=identity,$
(2) $D'=V(C')\subset D,$
(3) $\mu(C')\geq\mu(C)n^{-2L},$
(4) $\mu(D')\geq n^{-4L}e^{-16\omega(\phi)-8L\omega_0(\phi)}\mu(D)$
(5) $|\log(d\mu V/d\mu)(x)-\log(d\mu V/d\mu)(x')|<8\omega(\phi)$ *for all* $x,x'\in C'$.

*Proof.* For every $(i,j)\in\{0,\ldots n-1\}^2$ we denote by

$$S(i,j)\subset(\mathbb{Z}/n\mathbb{Z})^{L+1}$$

the set of allowed strings of the form $s=is_1\ldots s_{L-1}j$, and we write $c=i_{-M}\ldots i_M$ and $d=j_{-M}\ldots j_M$ for the strings defining the cylinders $C$ and $D$. If $(i,j)\in(\mathbb{Z}/n\mathbb{Z})^2$, $s\in S(i,i_{-M})$, $s'\in S(i_M,j)$, $t\in S(i,j_{-M})$ and $t'\in S(j_M,j)$, we consider the concatenations $scs',tdt'$ of length $2L+2M+1$ and denote by $[scs']_{-M-L}$ and $[tdt']_{-M-L}$ the corresponding cylinder sets. For every $x \in [scs']_{-M-L}$ we define $x' \in [tdt']_{-M-L}$ by

$$x'_k = \begin{cases} x_k & \text{if } |k| \geq M+L, \\ t_l & \text{if } k = -L-M+l, \, l = 0, \ldots, L, \\ j_k & \text{if } |k| \geq M, \\ t'_l & \text{if } k = M+l, \, l = 0, \ldots, L, \end{cases}$$

and observe that the map $x \to x'$ is a homeomorphism of $[scs']_{-L-M}$ onto $[tdt']_{-L-M}$. As the cylinders $C, D$ are disjoint this allows us to define an element $V_{s,s'}^{t,t'} \in [R_A]$ by setting

$$V_{s,s'}^{t,t'}(x) = \begin{cases} x & \text{if } x \in \Sigma_A \setminus ([scs']_{-L-M} \cup [tdt']_{-L-M}), \\ x' & \text{if } x \in [scs']_{-L-M}, \\ (V_{s,s'}^{t,t'})^{-1}(x) & \text{if } x \in [tdt']_{-L-M}. \end{cases}$$

From the definition of $\omega(\phi)$, comparing $\phi(\sigma^k x)$ to $\phi(\sigma^k x')$ (and $\phi(\sigma^k y)$ to $\phi(\sigma^k y')$) when $|k| \geq M+L$, and $\phi(\sigma^k x)$ to $\phi(\sigma^k y)$ (and $\phi(\sigma^k x')$ to $\phi(\sigma^k y')$) when $|k| < M+L$, it follows that

$$\max_{x,y \in [scs']} |\log \rho_\mu(V_{s,s'}^{t,t'}(x), x) - \log \rho_\mu(V_{s,s'}^{t,t'}(y), y)| < 8\omega(\phi). \tag{3.5}$$

In order to check the dependence of (3.5) on $s, s', t, t'$ we fix arbitrary points $\bar{x} \in C$, $\bar{y} \in D$ and calculate that

$$\left| \log \rho_\mu(V_{s,s'}^{t,t'}(x), x) - \sum_{k=-M-L}^{M+L} (\phi\sigma^k(\bar{y}) - \phi\sigma^k(\bar{x})) \right| < 8\omega(\phi) + 4L\omega_0(\phi) \tag{3.6}$$

for every $(i,j) \in (\mathbb{Z}/n\mathbb{Z})^2$, $s \in S(i, i_{-M})$, $s' \in S(i_M, j)$, $t \in S(i, j_{-M})$, $t' \in S(j_M, j)$ and $x \in [scs']_{-L-M}$. In other words, we have found a constant

$$\xi(D, C) = \sum_{k=-L-M}^{L+M} \phi\sigma^k(\bar{y}) - \phi\sigma^k(\bar{x})$$

such that the logarithm of the Radon-Nikodym derivative of each of the maps $V_{s,s'}^{t,t'}$ is within $8\omega(\phi) + 4L\omega_0(\phi)$ of this constant. In particular,

$$e^{-8\omega(\phi)-4L\omega_0(\phi)\xi(D,C)} < \frac{\mu([tdt']_{-L-M})}{\mu([scs']_{-L-M})} < e^{8\omega(\phi)+4L\omega_0(\phi)\xi(D,C)} \tag{3.7}$$

for every $(i,j) \in (\mathbb{Z}/n\mathbb{Z})^2$, $s \in S(i, i_{-M})$, $s' \in S(i_M, j)$, $t \in S(i, j_{-M})$ and $t' \in S(j_M, j)$. Since $\mu([scs']_{-L-M}) \geq \mu(C)n^{-2L}$ for at least one choice of $s, s'$ we conclude that, for any suitable choice of $t, t'$,

$$\mu(D) \geq \mu([tdt']_{-L-M}) \geq \mu([scs']_{-L-M})e^{-8\omega(\phi)-4L\omega_0(\phi)\xi(D,C)}$$
$$\geq \mu(C)n^{-2L}e^{-8\omega(\phi)-4L\omega_0(\phi)\xi(D,C)}.$$

Similarly we see that, for suitable choices of $s, s', t, t'$,

$$\mu(D) \leq \mu([tdt']_{-L-M})n^{2L} \leq \mu([scs']_{-L-M})n^{2L}e^{8\omega(\phi)+4L\omega_0(\phi)\xi(D,C)}$$
$$\leq \mu(C)n^{2L}e^{8\omega(\phi)+4L\omega_0(\phi)\xi(D,C)},$$

so that

$$n^{-2L}e^{-8\omega(\phi)-4L\omega_0(\phi)\xi(D,C)} \leq \frac{\mu(D)}{\mu(C)} \leq n^{2L}e^{8\omega(\phi)+4L\omega_0(\phi)\xi(D,C)}. \tag{3.8}$$

We choose $(i, j) \in (\mathbb{Z}/n\mathbb{Z})^2$ and $s \in S(i, i_{-M}), s' \in S(i_M, j)$ such that (3) holds for $C' = [scs']_{-L-M}$, choose $t \in S(i, j_{-M})$, $t' \in S(j_M, j)$, and set $V = V_{s,s'}^{t,t'}$.

Conditions (1)–(2) follow from the definition of $V$ with

$$D' = [tdt']_{-L-M},$$

(4) is a consequence of (3.7)–(3.8), and (5) is (3.5).     □

*Completion of the proof of Theorem 3.1.* Every probability measure in $M_1^\phi(\Sigma_A)$ has its Radon-Nikodym derivative under $R_A$ given by (3.4). In particular, by the chain rule any two distinct $R_A$-ergodic elements of $M_1^\phi(\Sigma_A)$ have to be mutually singular. Thus if $M_1^\phi(\Sigma_A)$ contains more than one probability measure, then it must therefore also contain a measure which is not $R_A$-ergodic. Choose a Borel set $B = R_A(B) \subset \Sigma_A$ with $0 < \mu(B) < 1$ and let

$$\varepsilon = 1/(100n^{4L}e^{24\omega(\phi)+8L\omega_0(\phi)}),$$

where $L \geq 0$ is chosen so that every entry of $A^L$ is positive.

Since $\mu(B)\mu(\Sigma_A \smallsetminus B) > 0$, there exist an integer $M \geq 0$ and cylinder sets $C = [i_{-M}, \ldots, i_M]_{-M}, D = [j_{-M}, \ldots, j_M]_{-M}$ in $\Sigma_A$ such that

$$\mu(B \cap C) > (1 - \varepsilon)\mu(C), \ \ \mu((\Sigma_A \smallsetminus B) \cap D) > (1 - \varepsilon)\mu(D);$$

we apply Lemma 3.2 to find cylinder sets $C' \subset C$ and $D' \subset D$, satisfying the conditions (1)–(5) of that lemma. According to (3),

$$\mu(C' \cap B) > \mu(C') - \varepsilon\mu(C) \geq \mu(C')(1 - \varepsilon n^{2L}),$$

and (5) guarantees that

$$
\begin{aligned}
\mu(V(C' \cap B)) &> \mu(D') e^{-8\omega(\phi)} (1 - \varepsilon n^{2L}) \\
&\geq \mu(D) n^{-4L} e^{-24\omega(\phi) - 8L\omega_0(\phi)} (1 - \varepsilon n^{2L}) \\
&= \mu(D) \frac{99}{100 n^{4L} e^{24\omega(\phi) + 8L\omega_0(\phi)}} \\
&> \varepsilon \mu(D) > \mu(B \cap D),
\end{aligned}
$$

which shows that $B$ cannot be invariant under $V$. This contradiction proves that $M_1^\phi(\Sigma_A)$ consists of a single $R_A$-ergodic measure, as claimed.

Finally, $M_1^\phi(\Sigma_A)$ is clearly shift-invariant, so the unique measure $\mu_\phi$ it contains is shift-invariant. Ergodicity of $\mu_\phi$ under $R_A$ is equivalent to triviality of the two-sided tail field, which implies triviality of the two one-sided tail fields, which in turn is equivalent to the $K$ property. □

For the remainder of this section we fix a map $\phi \colon \Sigma_A \longmapsto \mathbb{R}$ with summable variation and write $\mu_\phi$ for the unique Gibbs measure of $\phi$. According to Theorem 3.1, $\mu_\phi$ is ergodic under the Gibbs equivalence relation $R_A$; we shall now prove that $\mu_\phi$ is also ergodic under a class of shift-invariant subrelations of $R_A$ which was discussed in [48].

Let $\mathcal{G}$ be a countable, discrete group with finite conjugacy classes (i.e. with $|\{hgh^{-1} : h \in \mathcal{G}\}| < \infty$ for every $g \in \mathcal{G}$), and let $\psi \colon \Sigma_A \longmapsto \mathcal{G}$ be a continuous map.

Following [48] we set, for every $(x, x') \in R_A$,

$$
\begin{aligned}
\mathcal{J}_+^\psi(x, x') = \lim_{K \to \infty} \ &\psi(x) \cdots \psi(\sigma^K(x)) \\
&\cdot (\psi(x') \cdots \psi(\sigma^K(x')))^{-1}, \\
\mathcal{J}_-^\psi(x, x') = \lim_{K \to \infty} \ &(\psi(\sigma^{-K}(x)) \cdots \psi(\sigma^{-1}(x)))^{-1} \\
&\cdot \psi(\sigma^{-K}(x')) \cdots \psi(\sigma^{-1}(x'))
\end{aligned}
\tag{3.9}
$$

(since $\psi(\sigma^n(x)) = \psi(\sigma^n x')$ for $n$ outside a finite interval, these limits exist). The maps $\mathcal{J}_\pm^\psi \colon R_A \longmapsto \mathcal{G}$ are *cocycles* of $R_A$, i.e. satisfy that

$$
\mathcal{J}_\pm^\psi(x, x') \mathcal{J}_\pm^\psi(x', x'') = \mathcal{J}_\pm^\psi(x, x'')
\tag{3.10}
$$

whenever $(x, x'), (x, x'') \in R_A$. If $h \colon \Sigma_A \longmapsto \mathcal{G}$ is a continuous map,

$$
\psi'(x) = h(x) \psi(x) h(\sigma(x))^{-1},
\tag{3.11}
$$

and if $\mathcal{J}_\pm^{\psi'} \colon R_A \longmapsto \mathcal{G}$ are the cocycles defined as in (3.9), but with $\psi'$ replacing $\psi$, then

$$
\mathcal{J}_\pm^{\psi'}(x, x') = h(x) \mathcal{J}_\pm^\psi(x, x') h(x')^{-1}
\tag{3.12}
$$

for every $(x, x') \in R_A$ (this means that $\mathcal{J}_{\pm}^{\psi}$ and $\mathcal{J}_{\pm}^{\psi'}$ are *cohomologous* cocycles of $R_A$). The set

$$S_A^{\psi} = \{(x, x') \in R_A : \mathcal{J}_+^{\psi}(x, x') = \mathcal{J}_-^{\psi}(x, x')\} \qquad (3.13)$$

is a subrelation of $R_A$, and (3.12) shows that $S_A^{\psi}$ is unaffected if $\psi$ is replaced by $\psi'$ in (3.13). We shall prove the following theorem.

**Theorem 3.3.** *Let $\phi \colon \Sigma_A \longmapsto \mathbb{R}$ be a function with summable variation, and let $\mathcal{G}$ be a countable, discrete group with finite conjugacy classes. For every continuous map $\psi \colon \Sigma_A \longmapsto \mathcal{G}$ the Gibbs measure $\mu_\phi$ is ergodic under the equivalence relation $S_A^{\psi}$ defined in (3.13).*

The proof of Theorem 3.3 requires some preparation. Since $\mathcal{G}$ is discrete there exist integers $K, K' \geq 0$ such that $\psi(x)$ only depends on the coordinates $x_{-K}, \ldots, x_{K'}$ of every $x \in \Sigma_A$, and by applying a standard recoding argument (going to a higher-block presentation of $\Sigma_A$) we may assume without loss in generality that $K = 0, K' = 1$. We fix $L \geq 1$ such that every entry of $A^L$ is positive and find allowed strings $s^{(i)} = 0 s_1^{(i)} \ldots s_{L-1}^{(i)} i$ for every $i \in \{0, \ldots, n-1\}$. Define a continuous map $\bar{h} \colon \Sigma_A \longmapsto \mathcal{G}$ by setting

$$\bar{h}(x) = \psi(0, s_1^{(x_0)}) \cdots \psi(s_{L-1}^{(x_0)}, x_0) \qquad (3.14)$$

for every $x = (x_k) \in \Sigma_A$, put

$$\bar{\psi}(x) = \bar{h}(x)\psi(x)\bar{h}(\sigma(x))^{-1} \qquad (3.15)$$

for every $x \in \Sigma_A$, and denote by $\mathcal{J}_{\pm}^{\bar{\psi}}$ the cocycles defined by (3.9) with $\bar{\psi}$ replacing $\psi$.

**Lemma 3.4.** *Let $\Delta_{\pm} = \{\mathcal{J}_{\pm}^{\bar{\psi}}(x, x') : (x, x') \in R_A\}$.*

(1) *For every $(x, x') \in R_A$ there exist allowed strings $s = s_0 \ldots s_m$, $s' = s_0' \ldots s_m'$, $t = t_0 \ldots t_{m'}$, $t' = t_0' \ldots t_{m'}'$ such that $s_0 = s_m = s_0' = s_m' = t_0 = t_{m'} = t_0' = t_{m'}' = 0$ and*

$$\mathcal{J}_+^{\bar{\psi}}(x, x') = \psi(s_0, s_1) \cdots \psi(s_{m-1}, s_m)(\psi(s_0', s_1') \cdots \psi(s_{m-1}', s_m'))^{-1},$$

$$\mathcal{J}_-^{\bar{\psi}}(x, x') = (\psi(t_0, t_1) \cdots \psi(t_{m-1}, t_m))^{-1}\psi(t_0', t_1') \cdots \psi(t_{m-1}', t_m');$$

(2) *For every $m \geq 1$ and all allowed strings $s = s_0 \ldots s_m$, $s' = s_0' \ldots s_m'$ with $s_0 = s_m = s_0' = s_m' = 0$,*

$$\psi(s_0, s_1) \cdots \psi(s_{m-1}, s_m)(\psi(s_0', s_1') \cdots \psi(s_{m-1}', s_m'))^{-1} \in \Delta_+,$$

$$(\psi(s_0, s_1) \cdots \psi(s_{m-1}, s_m))^{-1}\psi(s_0', s_1') \cdots \psi(s_{m-1}', s_m') \in \Delta_-;$$

(3) *$\Delta_{\pm}$ is a subgroup of $\mathcal{G}$;*

(4) *For every $g_\pm \in \Delta_\pm$ and every $B \in \mathcal{B}_{\Sigma_A}$ with $\mu_\phi(B) > 0$ there exist elements $(x, x'), (y, y') \in R_A \cap (B \times B)$ with $x_k = x'_k$ for all $k \leq 0$, $y_k = y'_k$ for all $k \geq 0$, and*

$$g_+ = \mathcal{J}_+^{\bar\psi}(x, x'), \quad g_- = \mathcal{J}_-^{\bar\psi}(y, y').$$

*In fact we can accomplish this with $x = y$, $x' = y'$.*

*Proof.* We shall prove the conditions (1)–(4) for $\Delta_+$ and $\mathcal{J}_+^{\bar\psi}$; the proofs for $\Delta_-$ and $\mathcal{J}_-^{\bar\psi}$ are completely analogous and will be omitted.

Condition (1) follows from (3.12): if $(x, x') \in R_A$ satisfy that $x_k = x'_k$ for $|k| \geq K - L$, say, then

$$\begin{aligned}
\mathcal{J}_+^{\bar\psi}(x, x') &= \bar\psi(x_0, x_1) \cdots \bar\psi(x_{K-1}, x_K) \bar\psi(x'_{K-1}, x'_K)^{-1} \cdots \bar\psi(x'_0 x'_1)^{-1} \\
&= \psi(0, s_1^{(x_0)}) \cdots \psi(s_{L-1}^{(x_0)}, x_0) \psi(x_0, x_1) \cdots \psi(x_{K-1}, x_K) \\
&\quad \cdot \psi(x'_{K-1}, x'_K)^{-1} \cdots \psi(x'_0 x'_1)^{-1} \psi(s_{L-1}^{(x'_0)}, x'_0)^{-1} \cdots \psi(0, s_1^{(x'_0)})^{-1}.
\end{aligned}$$

Choosing a string $B$ such that $x_K B 0$ is allowed and setting

$$s = 0 s_1^{(x_0)} \ldots s_{L-1}^{(x_0)} x_0 \ldots x_K B 0, \quad s' = 0 s_1^{(x'_0)} \ldots s_{L-1}^{(x'_0)} x'_0 \ldots x'_K B 0$$

shows that $\mathcal{J}_+^{\bar\psi}(x, x')$ is of the required form.

In order to prove (4) and in the process (2) we assume that $s = s_0 \ldots s_m$, $s' = s'_0 \ldots s'_m$ are allowed strings strings with $s_0 = s'_0 = s_m = s'_m = 0$, set

$$g = \psi(s_0, s_1) \cdots \psi(s_{m-1}, s_m)(\psi(s'_0, s'_1) \cdots \psi(s'_{m-1}, s'_m))^{-1}$$

and denote by

$$\mathcal{C}(g) = \{h \in \mathcal{G} : hgh^{-1} = g\}$$

the commutant of $g$. As $\mathcal{G}$ has finite conjugacy classes, the quotient space $\mathcal{H} = \mathcal{G}/\mathcal{C}(g)$ is finite, and we define a homeomorphism $\bar\sigma : \Sigma_A \times \mathcal{H} \longmapsto \Sigma_A \times \mathcal{H}$ by setting

$$\bar\sigma(x, h\mathcal{C}(g)) = (\sigma(x), \bar\psi(x)^{-1} h\mathcal{C}(g))$$

for every $x \in \Sigma_A$ and $h \in \mathcal{G}$. We denote by $\nu$ the normalized counting measure on $\mathcal{H}$ and set $\bar\mu = \mu_\phi \times \nu$. If $B \in \mathcal{B}_{\Sigma_A}$ with $\mu_\phi(B) > 0$, then the mean ergodic theorem, applied to the indicator function of $B' = B \times \{\mathcal{C}(g)\} \subset \Sigma_A \times \mathcal{H}$, yields that

$$\lim_{l \to \infty} \frac{1}{l} \sum_{k=0}^{l-1} 1_{B'} \cdot \bar\sigma^k = E_{\bar\mu}(1_{B'} | \mathcal{B}_{\Sigma_A \times \mathcal{H}}^{\bar\sigma}) \quad \text{in } L^2(\bar\mu),$$

where $\mathcal{B}_{\Sigma_A \times \mathcal{H}}^{\bar{\sigma}}$ denotes the family of $\bar{\sigma}$-invariant Borel sets in $\Sigma_A \times \mathcal{H}$. Hence

$$\lim_{l \to \infty} \frac{1}{l} \sum_{k=0}^{l-1} \mu_\phi(B \cap \sigma^{-k}(B) \cap \{x \in \Sigma_A : \bar{\psi}(x) \cdots \bar{\psi}(\sigma^{k-1}(x)) \in \mathcal{C}(g)\})$$

$$= |\mathcal{H}| \cdot \lim_{l \to \infty} \frac{1}{l} \sum_{k=0}^{l-1} \bar{\mu}(B' \cap \bar{\sigma}^{-k}(B'))$$

$$\geq |\mathcal{H}| \cdot \bar{\mu}(B')^2 = \mu_\phi(B)^2 / |\mathcal{H}|.$$

In particular, there exist infinitely many $l \geq 0$ with

$$\mu_\phi(B \cap \sigma^{-l}(B) \cap \{x \in \Sigma_A : \bar{\psi}(x) \cdots \bar{\psi}(\sigma^{l-1}(x)) \in \mathcal{C}(g)\})$$
$$\geq \mu_\phi(B)^2 / 2|\mathcal{H}|. \tag{3.16}$$

Put $C = [s]_0$, $D = [s']_0$ and define $W \in [R_A]$ by

$$Wx = \begin{cases} (\ldots, x_{-1}, s'_0, \ldots, s'_m, x_{m+1}, \ldots) & \text{if } x = (x_k) \in C, \\ (\ldots, x_{-1}, s_0, \ldots, s_m, x_{m+1}, \ldots) & \text{if } x \in D, \\ x & \text{if } x \in \Sigma_A \smallsetminus (C \cup D). \end{cases}$$

Suppose that $B \in \mathcal{B}_{\Sigma_A}$, $B \subset [0]_0$ and $\mu_\phi(B) > 0$, and choose an $r \geq 1$ and an allowed string $t = t_0 \ldots t_{rm}$ such that $t_{rm} = 0$ and the set $\bar{B} = B \cap [t]_0$ has positive measure. We denote by $\bar{s} = \bar{s}_0 \ldots \bar{s}_{rm}$ the $r$-fold concatenation of $s$ and define $W' \in [R_A]$ by

$$W'x = \begin{cases} (\ldots, x_{-1}, \bar{s}_0, \ldots, \bar{s}_{rm}, x_{rm+1}, \ldots) & \text{if } x = (x_k) \in [t]_0, \\ (\ldots, x_{-1}, t_0, \ldots, t_{rm}, x_{rm+1}, \ldots) & \text{if } x \in [\bar{s}]_0, \\ x & \text{if } x \in \Sigma_A \smallsetminus ([t]_0 \cup [\bar{s}]_0). \end{cases}$$

Then, since

$$\lim_{l \to \infty} \mu_\phi(\bar{B} \cap \sigma^{-l} W \sigma^l(\bar{B})) = \mu_\phi(\bar{B}),$$

$$\lim_{l \to \infty} \mu_\phi(\bar{B} \cap \sigma^{-l}(\bar{B}) \cap \sigma^{-l} W \sigma^l(\bar{B}))$$
$$= \lim_{l \to \infty} \mu_\phi(\bar{B} \cap \sigma^{-l}(\bar{B}) \cap \sigma^{-l} W' \sigma^l(\bar{B}))$$
$$= \lim_{l \to \infty} \mu_\phi(\bar{B} \cap \sigma^{-l}(\bar{B}) \cap \sigma^{-l} WW' \sigma^l(\bar{B})) = \mu(\bar{B})^2,$$

and (3.16) guarantees that

$$\mu_\phi(\bar{B} \cap \sigma^{-l}(\bar{B}) \cap \sigma^{-l} W' \sigma^l(\bar{B}) \cap \sigma^{-l} WW' \sigma^l(\bar{B})$$
$$\cap \{x \in \Sigma_A : \bar{\psi}(x) \cdots \bar{\psi}(\sigma^{l-1}(x)) \in \mathcal{C}(g)\}) > 0 \tag{3.17}$$

for infinitely many $l \geq 0$. In particular, if $x = (x_k)$ lies in the set occurring in (3.17), then the points

$$x = (\ldots, x_{-1}, x_0, \ldots, x_l, t_1, \ldots, t_{rm}, x_{rm+1}, \ldots),$$
$$x' = (\ldots, x_{-1}, x_0, \ldots, x_l, \bar{s}_1, \ldots, \bar{s}_{rm}, x_{rm+1}, \ldots),$$
$$x'' = (\ldots, x_{-1}, x_0, \ldots, x_l, s'_1, \ldots, s'_m, \bar{s}_{m+1}, \ldots, \bar{s}_{r_m}, x_{rm+1}, \ldots)$$

all lie in $\bar{B} \subset B$, and

$$h = \bar{\psi}(x) \cdots \bar{\psi}(\sigma^{l-1}(x)) = \bar{\psi}(x') \cdots \bar{\psi}(\sigma^{l-1}(x'))$$
$$= \bar{\psi}(x'') \cdots \bar{\psi}(\sigma^{l-1}(x'')) \in \mathcal{C}(g).$$

Furthermore, $(x', x'') \in R_A \cap (B \times B)$, and

$$\mathcal{J}_+^{\bar{\psi}}(x', x'') = hgh^{-1} = g,$$

which proves (4) for any Borel set $B \subset [0]_0$ with positive measure.

If $B \in \mathcal{B}_{\Sigma_A}$ is an arbitrary Borel set with positive measure then there exists an allowed string $t' = t'_{-p} \ldots t'_p$ with $t'_p = t'_{-p} = 0$ and $\mu_\phi(B \cap [t']_{-p}) > 0$, since $\mu_\phi$ is K and hence mixing of every order. Put $B' = B \cap [t']_{-t}$ and $B'' = \sigma^p(B') \subset [0]_p$, and apply the preceding part of this proof to find, for every $g \in \Delta_+$, an element $(y, y') \in R_A \cap (B'' \times B'')$ with $y_k = y'_k$ for all $k \leq 0$ and $\mathcal{J}_+^{\bar{\psi}}(y, y') = g$. From the definition of $\mathcal{J}_+^{\bar{\psi}}$ it is clear that there exists an $h \in \Delta_+$ with $\mathcal{J}_+^{\bar{\psi}}(\sigma^{-p}(y), \sigma^{-p}(y')) = h\mathcal{J}_+^{\bar{\psi}}(y, y')h^{-1}$ for every $(y, y') \in R_A \cap (B'' \times B'')$. Note that $h$ does not depend on $g$. We conclude that there exists, for every $g \in \Delta_+$, an element $(x, x') \in R_A \cap (B' \times B')$ with $\mathcal{J}_+^{\bar{\psi}}(x, x') = hgh^{-1}$, yielding what is claimed in (4).

For the proof of (3) we assume that $g, h \in \Delta_+$ and apply (1) to find allowed strings $t = t_0 \ldots t_{m'}$, $t' = t'_0 \ldots t'_{m'}$, such that

$$h = \psi(t_0, t_1) \cdots \psi(t_{m-1}, t_m)(\psi(t'_0, t'_1) \cdots \psi(t'_{m-1}, t'_m))^{-1}.$$

According to (4) there exists an element $(x, x') \in R_A$ such that $x_k = x'_k$ for all $k \leq 0$ and $g = \mathcal{J}_+^{\bar{\psi}}(x, x')$. We choose $m \geq 0$ such that $x_k = x'_k$ whenever $k \geq m$, put $C = [x_0, \ldots, x_m]_0$, and apply (3.16) to find a point $y = (y_k) \in C$ and an integer $l > m$ such that $\sigma^l(y) \in C$ and

$$g' = \bar{\psi}(y) \cdots \bar{\psi}(\sigma^{l-1}(y)) \in \mathcal{C}(h).$$

Put

$$y' = (\ldots, y_{-1}, x_0, \ldots x_m, y_{m+1}, \ldots, y_{l-1}, t_0, \ldots, t_{m'}, z_{l+m'+1}, \ldots),$$
$$y'' = (\ldots, y_{-1}, x'_0, \ldots x'_m, y_{m+1}, \ldots, y_{l-1}, t'_0, \ldots, t'_{m'}, z_{l+m'+1}, \ldots),$$

where the coordinates $z_{l+m'+1}, z_{l+m'+2}, \ldots$ are arbitrary (but, of course, allowed), and obtain that $(y', y'') \in R_A$ and $\mathcal{J}_+^{\bar{\psi}}(y', y'') = gh$. $\qquad \square$

**Lemma 3.5.** *Let $\mathcal{J} = (\mathcal{J}_+^\psi, \mathcal{J}_-^\psi) \colon R_A \longmapsto \Delta = \Delta_+ \times \Delta_-$ be the cocycle*

$$\mathcal{J}(x, x') = (\mathcal{J}_+^\psi(x, x'), \mathcal{J}_-^\psi(x, x')),$$

*and let $\overline{R}$ be the equivalence relation on $\bar{X} = \Sigma_A \times \Delta$ defined by*

$$\overline{R} = \{((x, g, h), (x', g', h')) : (x, x') \in R_A,$$
$$g = \mathcal{J}_+^\psi(x, x')g', \ h = \mathcal{J}_-^\psi(x, x')h'\}.$$

*If $\nu$ is the counting measure on $\Delta$ and $\lambda = \mu_\phi \times \nu$, then $\lambda$ is quasi-invariant and ergodic under $\overline{R}$.*

*Proof.* The quasi-invariance of $\lambda$ under $\overline{R}$ is obvious. For every $(g, h) \in \Delta$ we define a homeomorphism $T_{(g,h)} \colon \bar{X} \longmapsto \bar{X}$ by setting

$$T_{(g,h)}(x, g', h') = (x, g'g, h'h)$$

for every $(x, g', h') \in \bar{X}$. For every Borel set $B \subset \bar{X}$ with $\lambda(B) > 0$ and $(g, h) \in \Delta$ we have that

$$T_{(g,h)}(\overline{R}(B)) = \overline{R}(T_{(g,h)}(B)).$$

Lemma 3.4(4) is easily seen to imply that every $\overline{R}$-saturated set is invariant under the transformations $T_{(g,h)}$, $(g, h) \in \Delta$, and the ergodicity of $R_A$ guarantees that $\lambda(\bar{X} \smallsetminus \overline{R}(B)) = 0$ for every Borel set $B \subset \bar{X}$ with $\lambda(B) > 0$ (cf. e.g. [59, 60, 46]). $\qquad \square$

*Proof of Theorem 3.3.* We use the notation of Lemmas 3.4–3.5 and denote by $1_\mathcal{G}$ the identity element in $\mathcal{G}$. Define a map $\eta \colon \Delta \longmapsto \mathcal{G}$ by $\eta(g, h) = gh^{-1}$ and set $\Xi = \eta(\Delta)$. We denote by $\overline{R}'$ the equivalence relation on $\bar{X}' = \Sigma_A \times \Xi$ given by

$$\overline{R}' = \{((x, g), (x', g')) : (x, x') \in R, \ g = \mathcal{J}_+^\psi(x, x')g' \mathcal{J}_-^\psi(x, x')^{-1}\}$$

and denote by $\bar{\xi} \colon \bar{X} \longmapsto \bar{X}'$ the map $\bar{\xi}(x, g, h) = (x, gh^{-1})$. If $\nu'$ is the counting measure on $\Xi$ and $\lambda' = \mu \times \nu'$, then $\bar{\xi} \colon (\bar{X}, \bar{\mu}) \longmapsto (\bar{X}', \lambda')$ is a nonsingular map which sends $\overline{R}$-equivalence classes to $\overline{R}'$-equivalence classes, and Lemma 3.5 implies that $\lambda'$ is ergodic under $\overline{R}'$. In particular, the equivalence relation $\overline{R}' \cap ((\Sigma_A \times \{1_\mathcal{G}\}) \times (\Sigma_A \times \{1_\mathcal{G}\}))$ induced by $\overline{R}'$ on $\Sigma_A \times \{1_\mathcal{G}\}$ is ergodic. Since $((x, 1_\mathcal{G}), (x', 1_\mathcal{G})) \in \overline{R}'$ if and only if $(x, x') \in S_A^\psi = S_A^{\psi'}$ we have proved Theorem 3.3. $\qquad \square$

*Remark 3.6.* As a special case of the fact that the equivalence relation $S_A^\psi$ in (3.13) is unaffected if $\psi$ is replaced by a function of the form $\psi' = h\psi(h \cdot \sigma)^{-1}$ in (3.11) we obtain that $S_A^\psi$ is shift-invariant: $(\sigma(x), \sigma(x')) \in$ $S_A^\psi$ if and only if $(x, x') \in S_A^\psi$. These properties of $S_A^\psi$ are particularly transparent if $\mathcal{G}$ is abelian. Define a map $\mathcal{J}^\psi \colon R_A \longmapsto \mathcal{G}$ by

$$\mathcal{J}^\psi(x, x') = \prod_{k=-\infty}^{\infty} (\psi(\sigma^k(x)\psi(\sigma^k(x')^{-1}) = \mathcal{J}_+^\psi(x, x')\mathcal{J}_-^\psi(x, x')^{-1} \tag{3.18}$$

for every $(x, x') \in R_A$. Then $\mathcal{J}^\psi$ is a cocycle, i.e. satisfies the equation (3.10) with $\mathcal{J}^\psi$ replacing $\mathcal{J}_\pm^\psi$, and

$$S_A^\psi = \{(x, x') \in R_A : \mathcal{J}^\psi(x, x') = 0\}.$$

As $\mathcal{J}^\psi = \mathcal{J}^{\psi'}$ whenever $\psi, \psi' \colon \Sigma_A \longmapsto \mathcal{G}$ are related via the equation (3.11), we obtain once again that $S_A^\psi = S_A^{\psi'}$. However, the individual cocycles $\mathcal{J}_\pm^\psi$ are obviously changed if we replace $\psi$ by $\psi'$.

## 4. Gibbs measures and subrelations of Gibbs equivalence relations on one-sided shift spaces

As in the preceding section we assume that $A = (A(i, j), 0 \leq i, j \leq n - 1)$ is an irreducible, aperiodic transition matrix with entries in $\{0, 1\}$. We define the one-sided SFT $\Sigma_A^+$ and the Gibbs relation $R_A^+$ as in Example 2.2 and write $\sigma = \sigma_A^+$ for the one-sided shift

$$\sigma(x)_k = x_{k+1}$$

on $\Sigma_A^+$. If $\phi \colon \Sigma_A^+ \longmapsto \mathbb{R}$ is a continuous map we define $\omega_k(\phi)$, $k \geq 0$, as in the two-sided case and say that $\phi$ has *summable variation* if $\sum_{k \geq 0} \omega_k(\phi) < \infty$. A probability measure $\mu$ on $\Sigma_A^+$ is a Gibbs measure of a continuous map $\phi \colon \Sigma_A \longmapsto \mathbb{R}$ if $\mu$ is quasi-invariant under $R_A^+$ and

$$\log \rho_\mu(x, x') = \sum_{k \geq 0} (\phi\sigma^k(x) - \phi\sigma^k(x')) \tag{4.1}$$

for every $(x, x') \in R_A^+$ (cf. (3.4); in the one-sided case this definition makes sense even if $\phi$ does not have summable variation). The same argument as in Theorem 3.1 shows that the set $M_1^\phi(\Sigma_A^+) \subset M_1(\Sigma_A^+)$ of Gibbs measures of $\phi$ is nonempty, weak*-closed, and convex. Furthermore, since every function $\phi \colon \Sigma_A^+ \longmapsto \mathbb{R}$ with summable variation can be viewed as a function on $\Sigma_A$ which again has summable variation, $M_1^\phi(\Sigma_A^+)$ contains a single measure $\mu_\phi$ which is ergodic under the Gibbs relation $R_A^+$. In contrast to the two-sided case the measure $\mu_\phi$ need not be shift-invariant, although Theorem 3.1 guarantees that it is equivalent to a shift-invariant and exact probability measure $\bar\mu_\phi$ on $\Sigma_A^+$. For example if $P = (P(i, j), 0 \leq i, j \leq n - 1)$ is a stochastic matrix which is compatible with $A$ in the sense of Section 3, and if $\bar p$ is the unique probability vector satisfying $\bar{p}P = \bar{p}$ and $p(i) = 1/n$ for $i = 0, \dots, n-1$, then

$$\mu_{\log P}(C) = p(v_0) P(v_0, v_1) \cdots P(v_{l-1}, v_l) \tag{4.2}$$

for every cylinder set

$$C = [v_0, \dots, v_l]_m = \{x = (x_k) \in \Sigma_A^+ : x_{m+i} = v_i \text{ for } i = 0, \dots, l\},$$

whereas the shift-invariant measure $\bar{\mu}_{\log P}$, given by

$$\bar{\mu}_{\log P}(C) = \bar{p}(v_0) P(v_0, v_1) \cdots P(v_{k-1}, x_k), \tag{4.3}$$

is the Gibbs measure of the function $\phi(x) = \log P(x_0, x_1) + \log \bar{p}(x_0) - \log \bar{p}(x_1)$ (cf. [41]). [The reason for the difference between the one- and two-sided cases is the following. If $\phi$ is changed to a cohomologous function, the Radon-Nikodym derivative of the two-sided measure under the finite coordinate changes is unaffected, but this is no longer true in the one-sided case. The uniqueness in the one-sided case is proved exactly as in the two-sided case; however, since the Radon-Nikodym derivative changes on the one-sided shift if $\phi$ is changed by a cohomologous function, the one-sided measure changes as well. The condition of shift-invariance in the one-sided case singles out a particular element in the cohomology class of $\phi$.] For easier reference we summarize this discussion in a theorem.

**Theorem 4.1.** *Let $A = (A(i, j), 0 \leq i, j \leq n - 1)$ be an irreducible, aperiodic 0-1-matrix, $\Sigma_A^+$ the associated one-sided shift of finite type, and $R_A^+$ the Gibbs relation of $\Sigma_A^+$. If $\phi \colon \Sigma_A^+ \longmapsto \mathbb{R}$ is a function with summable variation, then there exists a unique Gibbs measure $\mu_\phi$ of $\phi$ which is ergodic under $R_A$. Furthermore, $\mu_\phi$ is equivalent to a probability measure $\bar{\mu}_\phi$ which is invariant, ergodic and exact under the shift $\sigma$ on $\Sigma_A^+$.*

*Remark 4.2.* According to [58], the measure $\bar{\mu}_\phi$ is actually *Bernoulli* for the shift $\sigma$.

We encounter a similar phenomenon when dealing with cocycle-generated subrelations of the Gibbs relation $R_A^+$ on the one-sided SFT $\Sigma_A^+$. Fix a function $\phi \colon \Sigma_A^+ \longmapsto \mathbb{R}$ with summable variation, denote by $\mu_\phi$ the Gibbs measure of $\phi$, and consider a continuous map $\psi \colon \Sigma_A^+ \longmapsto \mathcal{G}$, where $\mathcal{G}$ is a countable, discrete group with finite conjugacy classes. We define a cocycle $\mathcal{J}_+^\psi$ by the first equation in (3.9) and put

$$S_A^{\psi+} = \{(x, x') \in R_A^+ : \mathcal{J}_+^\psi(x, x') = 1_{\mathcal{G}}\}, \tag{4.4}$$

where $1_{\mathcal{G}}$ is the identity element in $\mathcal{G}$. In contrast to the situation described in Remark 3.6, the equivalence relation $S_A^{\psi+}$ is changed if $\psi$ is replaced by a function $\psi' = h\psi(h \cdot \sigma)^{-1}$, where $h\colon \Sigma_A^+ \longmapsto \mathcal{G}$ is continuous. In particular, if $\bar{h}\colon \Sigma_A^+ \longmapsto \mathcal{G}$ is the function defined in (3.14) and $\bar{\psi}\colon \Sigma_A^+ \longmapsto \mathcal{G}$ is given by (3.15), then essentially the same proof as that of Theorem 3.3 yields the following result.

**Theorem 4.3.** *Let* $\phi\colon \Sigma_A^+ \longmapsto \mathbb{R}$ *be a function with summable variation, and let* $\mathcal{G}$ *be a countable, discrete group with finite conjugacy classes. For every continuous map* $\psi\colon \Sigma_A^+ \longmapsto \mathcal{G}$ *there exists a continuous map* $h\colon \Sigma_A^+ \longmapsto \mathcal{G}$ *such that the Gibbs measure* $\mu_\phi$ *is ergodic under the equivalence relation* $S_A^{\psi'+}$ *defined by* (4.4) *with*

$$\psi' = h\psi(h \cdot \sigma)^{-1}$$

*replacing* $\psi$.

*The measure* $\mu_\phi$ *is ergodic under the original relation* $S_A^{\psi+}$ *if and only if* $S_A^{\psi+}$ *is topologically transitive.*

*Proof.* The proof of Theorem 3.3 shows that the first assertion is a consequence of Lemma 3.5. The topological transitivity of $S_A^{\psi+}$ is obviously necessary for the ergodicity of $\mu_\phi$ under $S_A^{\psi+}$, since every open subset of $\Sigma_A$ has positive $\mu_\phi$-measure. In order to see that topological transitivity is also sufficient for ergodicity we observe that there exists a finite partition $\mathcal{O}_1, \ldots, \mathcal{O}_m$ of $\Sigma_A^+$ into closed and open subsets on each of which the continuous map $h\colon \Sigma_A^+ \longmapsto \mathcal{G}$ is constant, and that $S_A^{\psi'+} \cap (\mathcal{O}_i \times \mathcal{O}_i) = S_A^{\psi+} \cap (\mathcal{O}_i \times \mathcal{O}_i)$ for $i = 1, \ldots, m$ (cf. (3.12)). The ergodicity of $\mu_\phi$ under $S_A^{\psi'+}$ guarantees that the restriction of $\mu$ to $\mathcal{O}_i$ is ergodic under $S_A^{\psi+} \cap (\mathcal{O}_i \times \mathcal{O}_i)$, and the topological transitivity of $S_A^{\psi+}$ implies that $S_A^{\psi+}(\mathcal{O}_i)$ is a dense, open subset of $\Sigma_A^+$ for every $i = 1, \ldots, m$ (note that there exists a countable group $G$ of homeomorphisms of $\Sigma_A^+$ with $S_A^{\psi+} = Gx = \{gx : g \in G\}$ for every $x \in \Sigma_A$). In particular, $S_A^{\psi+}(\mathcal{O}_1)$ meets each $\mathcal{O}_i$ in a set of positive measure, so that $S_A^{\psi+}$ is ergodic. $\qquad\square$

*Remark 4.4.* Assume for simplicity that $\psi$ is a function of the two variables $(x_0, x_1)$, and call two symbols $a, b \in \{0, \ldots, n-1\}$ *equivalent* if there exist allowed strings $s, t$ of equal length and a symbol $i \in \mathbb{Z}/n\mathbb{Z}$ such that the strings $asi$ and $bti$ are allowed and are permutations of each other. Then it is easy to see that $S_A^{\psi+}$ is topologically transitive if and only if all symbols in $\mathbb{Z}/n\mathbb{Z}$ are equivalent (cf. [8, 9]). Thus ergodicity of $\mu_\phi$ under $S_A^{\psi+}$ depends only on $A$, and not on $\phi$ or $\psi$.

*Remark 4.5.* The results in sections 3–4 are largely unaffected if we drop the assumption of aperiodicity of $A$. If $A$ is irreducible, but no longer aperiodic, then the alphabet $\{0, \ldots, n-1\}$ decomposes into periodic components $C_0, \ldots, C_{m-1}$, $m \geq 2$, such that $x_{i+k} \in C_{j+k \pmod m}$ whenever $x = (x_k) \in \Sigma_A$, $x_i \in C_j$, and $k \in \mathbb{Z}$. If $X_j = \{x \in \Sigma_A : x_0 \in C_j\}$, $j = 0, \ldots, m-1$, then each $X_j$ is invariant under $\sigma^m$, and a standard recoding argument allows us to regard $X_j$ as an irreducible and aperiodic SFT with regards to the shift $\sigma^m$. By applying Theorem 3.1 to each $X_j$ we see that there exists, for every function $\phi \colon \Sigma_A \longmapsto \mathbb{R}$ with summable variation, and for every $j = 0, \ldots, m-1$, a unique Gibbs probability measure $\mu_\phi^{(j)}$ for $\phi$ on $X_j$. In particular, there is no longer a unique probability measure on $\Sigma_A$ satisfying (3.4).

The orbit-average of the Gibbs measure $\mu_\phi^{(0)}$, say, under $\sigma$ is the unique *shift-invariant* Gibbs measure $\mu_\phi$ of $\phi$ on $\Sigma_A$ (it coincides, of course, with the orbit-averages of the other measures $\mu_\phi^{(j)}$, $j = 1, \ldots, m-1$). If we define the Gibbs relation $R_A$ exactly as in the aperiodic case, then $R_A(X_j) = X_j$, and the analogue of Theorem 3.3 holds for each of the measures $\mu_\phi^{(j)}$ on $X_j$. However, $\mu_\phi$ will no longer be $S_A^\psi$-invariant for any $\psi$. Similarly one obtains the corresponding statements on one-sided SFT's.

## 5. Examples

Throughout this section we assume that $A = (A(i,j), 0 \leq i, j \leq n-1)$ is an irreducible and aperiodic 0-1-matrix, write $\Sigma_A \subset (\mathbb{Z}/n\mathbb{Z})^{\mathbb{Z}}$ and $\Sigma_A^+ \subset (\mathbb{Z}/n\mathbb{Z})^{\mathbb{N}}$ for the two- and one-sided SFT's defined by $A$, and denote by $R_A$ and $R_A^+$ the Gibbs relations on $\Sigma_A$ and $\Sigma_A^+$.

*Example 5.1. Invariant measures for equivalence relations.* Let $\phi \colon \Sigma_A \longmapsto \mathbb{R}$ be a continuous map which takes only finitely many values, and let $\mathcal{G} \subset \mathbb{R}$ be the countable subgroup generated by the values of $\phi$, furnished with the discrete topology. We put $\psi = \phi \colon \Sigma_A \longmapsto \mathcal{G}$ and apply Theorem 3.3 to obtain that the Gibbs measure $\mu_\phi$ is ergodic, and obviously invariant, under the equivalence relation $S_A^\psi \subset R_A$. Note that the full group $[S_A^\psi]$ satisfies that

$$[S_A^\psi] = \{W \in [R_A] : W \text{ preserves } \mu_\phi\}.$$

One special case of this construction is obtained by taking a stochastic matrix $P = (P(i,j), 0 \leq i, j \leq n-1)$ compatible with $A$ and setting $\phi(x) = \log P(x_0, x_1)$ for every $x \in \Sigma_A$. As noted in the discussion preceding Theorem 3.1, in this case $\mu_\phi$ is the shift-invariant Markov probability measure determined by $P$. (So every Markov measure is a $\mu_\phi$ in this way.)

Another special case arises from $\phi = \psi = $ constant. Then $\mu_\phi$ is the unique measure of maximal entropy on $\Sigma_A$. In the one-sided case, it is equivalent to the unique invariant measure for the stationary adic on $\Sigma_A^+$.

*Example 5.2. Failure of ergodicity.* In order to illustrate Theorem 4.1 and Remark 4.4, consider the one-sided two-shift $\Sigma_A^+ = (\mathbb{Z}/2\mathbb{Z})^{\mathbb{N}}$ and the map $\psi \colon \Sigma_A^+ \longmapsto \mathbb{Z}$ defined by

$$\psi(x) = x_0 - x_1$$

for every $x = (x_k) \in \Sigma_A^+$. Then

$$\mathcal{J}_+^\psi(x, x') = x_0 - x_0'$$

for every $(x, x') \in R_A^+$ (cf. (3.9)), and the relation $S_A^{\psi+}$ is obviously not topologically transitive. In particular, if $\phi \colon \Sigma_A^+ \longmapsto \mathbb{R}$ is a function with summable variation, then the Gibbs measure $\mu_\phi$ is nonergodic under $S_A^{\psi+}$.

*Example 5.3. Ergodicity of Gibbs measures on one-sided SFT's under finite permutations.* Suppose that $S_A^+ \subset R_A^+$ is the equivalence relation in which two points $x = (x_k)$, $x' = (x_k')$ in $\Sigma_A^+$ are equivalent if and only if the coordinates $x_k$, $k \geq 0$, and $x_k'$, $k \geq 0$, are finite permutations of each other (cf. Examples 2.1–2.2). If $\psi_0 \colon \Sigma_A^+ \longmapsto \mathbb{Z}^{\mathbb{Z}/n\mathbb{Z}}$ is the map defined by

$$\psi_0(x) = e^{(x_0)}$$

for every $x = (x_k) \in \Sigma_A^+$, where $e^{(0)} = (1, 0, \ldots, 0), \ldots, e^{(n-1)} = (0, \ldots, 0, 1)$ are the unit vectors in $\mathbb{Z}^{\mathbb{Z}/n\mathbb{Z}} \cong \mathbb{Z}^n$, then it is clear that $S_A^+ = S_A^{\psi_0+}$.

If

$$A = \left( \begin{smallmatrix} 1 & 1 \\ 1 & 1 \end{smallmatrix} \right),$$

then $S_A^+$ is obviously topologically transitive, since the strings 010 and 100 are allowed and permutations of each other (cf. Remark 4.4). However, if

$$A = \left( \begin{smallmatrix} 1 & 1 & 0 \\ 0 & 0 & 1 \\ 1 & 1 & 0 \end{smallmatrix} \right),$$

then $S_A^+$ is not topologically transitive. Indeed, if $s = s_0 s_1 \ldots s_m$ is any cycle (i.e. any allowed string with $s_0 = s_m$), then the numbers of 1's and 2's among the entries $s_0, \ldots, s_{m-1}$ must be equal. In particular, if the symbols 1 and 2 were equivalent in the sense of Remark 4.4, we could find allowed strings $s$ and $t$ of equal length and a symbol $i \in \{0, 1, 2\}$ such that $1si$ and $2ti$ are allowed and permutations of each other. We may obviously assume that $i = 1$; then the number of 1's in $1s1$ exceeds the number of 2's by 1, whereas the number of 1's and 2's in $2t1$ is equal. This shows that $1s1$ and $2t1$ cannot be permutations of each other, so that 1 and 2 are not equivalent and $S_A^{\psi+}$ is not topologically transitive (cf. Theorem 4.3 and Remark 4.4).

If $\phi: \Sigma_A^+ \longmapsto \mathbb{R}$ is a function with summable variation, then the Gibbs measure $\mu_\phi$ on $\Sigma_A^+$ is ergodic under $S_A^+$ only if $S_A^+$ is topologically transitive (Theorem 4.3). All these statements have, of course, obvious translations into properties of the adic transformations $T_A^*$ on $\Sigma_A^+$ and $T_{Y_A}$ on $Y_A$ (cf. Example 2.2).

*Example 5.4. Exchangeability and partial exchangeability.* The preceding example is a special case of a general problem: given some group of symmetries of a one- or two-sided SFT, what are the Gibbs measures which are invariant (or *symmetric*) under this group of symmetries? The most elementary examples of such symmetry groups are the group of all finite permutations of the coordinates $x_k$, $k \geq 0$, and the group of all finite permutations of the transitions $[x_k, x_{k+1}]$ of the points $x = (x_k)$ in $\Sigma_A$. We have met the first of these groups in the context of adic transformations in the Examples 2.1–2.2 and in Section 2, and our analysis of them depends on the fact that they generate the orbits (= equivalence classes) of certain subrelations of $R_A$ and $R_A^+$ defined by cocycles. In order to describe these cocycles we consider the group $\mathcal{G}_{l+1} = \mathbb{Z}^{\{0,\dots,n-1\}^{l+1}}$, denote by $e^{(i_0,\dots,i_l)} \in \mathbb{Z}^{(\mathbb{Z}/n\mathbb{Z})^{l+1}}$ the unit vector

$$e_{(j_0,\dots,j_l)}^{(i_0,\dots,i_l)} = \begin{cases} 1 & \text{if } (i_0,\dots,i_l) = (j_0,\dots,j_l), \\ 0 & \text{otherwise,} \end{cases}$$

and observe that the orbits of the group of all permutations of the $l$-fold transitions $[x_k,\dots,x_{k+l}]$, $k \geq 0$, of the points in $\Sigma_A$ and $\Sigma_A^+$ are the equivalence classes of the relations $S_A^{\psi_l}$ and $S_A^{\psi_l+}$, where $\psi_l$ is defined by

$$\psi_l(x) = e^{(x_0,\dots,x_l)} \in \mathcal{G}_{l+1}$$

for every $x$ in $\Sigma_A$ or $\Sigma_A^+$.

Theorem 3.3 shows, for a two-sided SFT $\Sigma_A$, the invariance and ergodicity of every fully supported Markov measure under the group of all finite permutations of the $l$-fold transitions described above. Conversely, Theorem 6.2 implies that every shift-invariant probability measure on $\Sigma_A$ which is invariant and ergodic under the group of finite permutations of the $l$-fold transitions must be an $l$-step Markov measure (where we are referring to the Bernoulli measures as zero-step Markov measures).

These results carry over to the one-sided case with the precaution described in Remarks 4.4 and 4.5 and Examples 5.1 and 5.3. The adic version of the statement of ergodicity and invariance of Markov measures under the groups of all finite permutations of the $l$-fold transitions is the Hewitt-Savage 0-1-Law as extended to Markov chains by Diaconis and Freedman [8].

*Example 5.5. Derangement-equivalent sequences.* In all the examples so far the group $\mathcal{G}$ has been abelian. An elementary example where $\mathcal{G}$ is nonabelian is obtained by setting $\mathcal{G}$ equal to the group $\mathfrak{S}_n$ of permutations of the symbols $\mathbb{Z}/n\mathbb{Z}$, and by putting

$$\psi(x) = (x_0 x_1)$$

for every $x = (x_k) \in \Sigma_A$, where $(ij) \in \mathfrak{S}_n$ is the transposition of $i$ and $j$. Since $\mathfrak{S}_n$ is finite, Theorem 3.3 implies that, for every function $\phi \colon \Sigma_A \longmapsto \mathbb{R}$ with summable variation, the Gibbs measure $\mu_\phi \in M_1(\Sigma_A)$ is ergodic under the relation $S = S_A^\psi$. In Examples 5.3 and 5.4, two sequences are equivalent if they eventually accumulate the same symbol counts (or perhaps the same transition counts). Here two sequences are equivalent if they eventually accumulate the same *derangement effects*, when their symbols act as permutations on a fixed set of $n$ letters.

The equivalence relation $S$ defined here is not comparable with the relation $S_A^{\psi_1}$ arising from the finite permutations of the transitions $[x_k, x_{k+1}]$ appearing in the Examples 5.3–5.4. However, since every finite Cartesian product of groups with finite conjugacy classes is again a group with finite conjugacy classes, Theorem 3.3 guarantees the following: if $\phi \colon \Sigma_A \longmapsto \mathbb{R}$ is a function with summable variation and $\psi^{(i)} \colon \Sigma_A \longmapsto \mathcal{G}^{(i)}$ are continuous maps taking values in countable, discrete groups with finite conjugacy classes $\mathcal{G}^{(i)}$, $i = 1, \ldots, l$, then the Gibbs measure $\mu_\phi$ is ergodic under the equivalence relation $\bigcap_{i=1}^l S_A^{\psi^{(i)}}$. In particular, $\mu_\phi$ is ergodic under $S_A^{\psi_1} \cap S$, the relation that keeps track of *both* accumulated symbol accounts and accumulated derangement effects.

As a concrete example, consider the case where $n = 4$ and

$$A = \begin{pmatrix} 1 & 1 & 1 & 1 \\ 1 & 1 & 1 & 1 \\ 1 & 1 & 1 & 1 \\ 1 & 1 & 1 & 1 \end{pmatrix},$$

so that $\Sigma_A = \Sigma_4$ is the full 4-shift. The points

$$\ldots x_{-1} 0123020 x_7 \ldots, \qquad \ldots x_{-1} 0230120 x_7 \ldots$$

are $S_A^{\psi_1}$-equivalent, but $S$-inequivalent, whereas the points

$$\ldots x_{-1} 010 x_3 \ldots, \qquad \ldots x_{-1} 000 x_3 \ldots$$

are $S$-equivalent, but $S_A^{\psi_1}$-inequivalent.

## 6. Symmetric Measures on SFT's

We continue to let $A = (A(i, j), 0 \leq i, j \leq n - 1)$ be an irreducible and aperiodic 0-1-matrix, $\Sigma_A$ and $\Sigma_A^+$ the associated two- and one-sided SFT's, and $S_A \subset R_A$ the equivalence relation generated by the finite coordinate permutations (cf. Example 2.2). As is well known, the only probability measure $\mu$ on $\Sigma_A$ which is invariant under the Gibbs relation $R_A$ is the measure of maximal entropy. However, if $R \subset R_A$ is a proper subrelation, the picture changes considerably. For example, if $\Sigma_A$ is a full shift, then de Finetti's theorem states that the (shift-invariant) Bernoulli measures are precisely the $S_A$-invariant measures. If $A$ is arbitrary (but still irreducible and aperiodic), and if $\phi \colon \Sigma_A \longmapsto \mathbb{R}$ is a function depending on the single coordinate $x_0$, then the equation (3.4) shows that the Gibbs measure $\mu_\phi \in M_1(\Sigma_A)$ is invariant under $S_A$, and the ergodicity of $\mu_\phi$ follows from Theorem 3.3. In Theorem 2.11 we saw that for the golden mean shift $\Sigma_A$, every $S_A$-invariant probability measure on $\Sigma_A^+$ was the Gibbs measure of a function of a single variable. In this section we extend the statement to arbitrary SFT's by showing that every shift- and $S_A$-invariant probability measure on $\Sigma_A$ is the Gibbs measure for a potential function that depends on just one coordinate. Very closely related results have been obtained previously by Georgii [14] and Diaconis and Freedman [8]. By an obvious extension to $m$-step SFT's, every shift- and $S_A^{\psi_m}$-invariant probability measure $\mu$ is the Gibbs measure of a function depending on the $m + 1$ variables $x_0, \dots, x_m$.

We write $\pi_\mathbb{N} \colon \Sigma_A \longmapsto \Sigma_A^+$ for the projection onto the nonnegative coordinates and observe that, in the notation of Example 2.2,

$$\bigcup_{i=0}^{n-1} S_A^+ \cap ([i]_0 \times [i]_0) \subset S = (\pi_\mathbb{N} \times \pi_\mathbb{N})(S_A) \subset R_A^+. \tag{6.1}$$

**Lemma 6.1.** *Let $\mu$ be a nonatomic, shift-invariant probability measure on $\Sigma_A$ which is invariant and ergodic under $S_A$. Then the projection $\mu^+ = \mu\pi_\mathbb{N}^{-1}$ of $\mu$ onto $\Sigma_A^+$ is quasi-invariant and ergodic under $S = (\pi_\mathbb{N} \times \pi_\mathbb{N})(S_A)$. Furthermore, the restriction of $\mu^+$ to the set $[i]_0 = \{x \in \Sigma_A^+ : x_0 = i\}$ is invariant and ergodic under the relation $S_A^+ \cap ([i]_0 \times [i]_0)$ for every $i \in \{0, \dots, n - 1\}$ with $\mu^+([i]_0) > 0$.*

*Proof.* The quasi-invariance and ergodicity of $\mu^+$ under $S$ are obvious from (6.1) and the definition of $\mu^+$. Suppose that $\mu^+([i]_0) > 0$, and that there exist disjoint Borel sets $B_1, B_2 \subset [i]_0$ with positive measure, each of which is saturated under $S_A^+ \cap ([i]_0 \times [i]_0)$. The ergodicity of $\mu^+$ under $S$ implies that there exists a homeomorphism $V \in [S]$ with $\mu^+(VB_1 \cap B_2) > 0$. In particular we can find a $j \in \{0, \dots, n - 1\}$ and allowed strings $s^{(1)} = is_1^{(1)} \ldots s_{m-1}^{(1)} j$, $s^{(2)} = is_1^{(2)} \ldots s_{m-1}^{(2)} j$, such that $V([s^{(1)}]_0) = [s^{(2)}]_0$ and the sets

$$C_l = B_l \cap [s^{(l)}]_0, \;\; l = 1, 2,$$

have positive measure. We put $s_0^{(l)} = i$, $s_m^{(l)} = j$ for $l = 1, 2$ and define, for every $p \geq 0$, $V^{(p)} \in [R_A]$ by

$$(V^{(p)} x)_k = \begin{cases} s_{k-p}^{(2)} & \text{if } \sigma^p(x) \in C_1 \text{ and } p \leq k \leq p + m, \\ s_{k-p}^{(1)} & \text{if } \sigma^p(x) \in C_2 \text{ and } p \leq k \leq p + m, \\ x_k & \text{if } k < p, \; k > p + m, \text{ or } \sigma^p(x) \notin (C_1 \cup C_2). \end{cases}$$

In other words, $V^{(p)}$ checks for each $x \in \Sigma_A^+$ whether either of the strings $s^{(0)}, s^{(1)}$ occurs at the positions $p, \ldots, p + m$; if one of these strings occurs, then it replaces it by the other one, and if neither occurs, then the point is left unchanged. The measure $\mu^+$ is quasi-invariant under each of the maps $V^{(p)}$, and the shift-invariance of $\mu^+$ guarantees that the sequence of Radon-Nikodym derivatives $(d\mu^+ V^{(p)} / d\mu^+, \; p \geq 0)$ is uniformly integrable. By approximating $C_i$ by closed and open sets and using uniform integrability we see that

$$\lim_{p \to \infty} \mu^+(C_l^{(p)}) = \mu^+(C_l)$$

for $l = 1, 2$, where

$$C_l^{(p)} = C_l \cap V^{(p)} C_l.$$

The ergodicity of $\mu^+$ under $S$ implies that $\mu^+$ is exact and hence mixing under the shift $\sigma$ on $\Sigma_A^+$, so that

$$\lim_{p \to \infty} \mu^+(C_l \cap [s^{(k)}]_p) = \mu^+(C_l) \mu^+([s^{(k)}]_0)$$

for every $k, l \in \{0, 1\}$. Hence

$$\lim_{p \to \infty} \mu^+(C_l^{(p)} \cap [s^{(k)}]_p) = \mu^+(C_l) \mu^+([s^{(k)}]_0)$$

for every $k, l \in \{0, 1\}$, and by setting

$$W^{(p)} x = \begin{cases} V^{(0)} V^{(p)} x & \text{for } x \in (C_1^{(p)} \cap [i^{(2)}]_p) \cup (C_2^{(p)} \cap [i^{(1)}]_p) \\ x & \text{otherwise} \end{cases}$$

for every $p > m$ we have constructed a sequence of maps in $[S_A^+]$ with

$$\lim_{p \to \infty} \mu^+(C_2 \cap W^{(p)} C_1) = \mu^+(C_2) \mu^+([s^{(2)}]_0) > 0.$$

This violates the invariance of $B_1$ and $B_2$ under $S_A^+ \cap ([i]_0 \times [i]_0)$, and we conclude that the restriction of $\mu^+$ to $[i]_0$ is ergodic under $S_A^+ \cap ([i]_0 \times [i]_0)$. $\qquad\square$

With this lemma at hand we can characterize the set of all shift-invariant measures on $\Sigma_A$ which are invariant and ergodic under $S_A$. Any 1-step (irreducible aperiodic) $S_A$-invariant SFT inside $\Sigma_A$ with a potential function on it depending on only a single coordinate determines such a measure, and we will show that they all arise in this way. In the following theorem we denote by $\Sigma_{A'} \subset \Sigma_A$ an irreducible and aperiodic subshift defined by a 0-1-matrix $A' = (A'(i,j), \, 0 \leq i, j \leq n-1)$ with $A'(i,j) \leq A(i,j)$ for every $(i,j) \in \{0, \ldots, n-1\}^2$. The matrix $A'$ may have some zero rows and columns; however, the irreducibility of $\Sigma_{A'}$ allows us to permute the alphabet $\{0, \ldots, n-1\}$ of $\Sigma_A$, if necessary, and to assume that there exists an integer $n' \in \{1, \ldots, n-1\}$ such that the first $n'$ rows and columns of $A$ are nonzero and the remaining ones zero, and to regard $\Sigma_{A'}$ as a (irreducible and aperiodic) subshift of $\{0, \ldots, n'-1\}^{\mathbb{Z}}$. Any probability measure $\mu$ on $\Sigma_A$ with $\mu(\Sigma_{A'}) = 1$ will also be regarded as a probability measure on the SFT $\Sigma_{A'} \subset \Sigma_A$, and vice versa.

**Theorem 6.2.** *Let $\mu$ be a nonatomic shift-invariant probability measure on $\Sigma_A$ which is invariant and ergodic under $S_A$. Then there exists an irreducible, aperiodic and $S_A$-invariant SFT $\Sigma_{A'} \subset \Sigma_A$ and a function $\phi \colon \Sigma_{A'} \to \mathbb{R}$ depending on the single coordinate $x_0$ of every point $x \in \Sigma_{A'}$ such that $\mu$ is the unique (Markov) Gibbs measure of $\phi$ on $\Sigma_{A'}$.*

*Conversely, if $\Sigma_{A'} \subset \Sigma_A$ is an irreducible, aperiodic and $S_A$-invariant SFT and $\phi \colon \Sigma_{A'} \longmapsto \mathbb{R}$ is a function of a single coordinate, then the (Markov) Gibbs measure $\mu_\phi$ is invariant and ergodic under $S_A$.*

*Proof.* First we construct an irreducible and aperiodic SFT $\Sigma_{A'} \subset \Sigma_A$ with $\mu(\Sigma_{A'}) = 1$, such that $\mu$ is the Gibbs measure of a function $\phi \colon \Sigma_{A'} \longmapsto \mathbb{R}$ of only two coordinates.

Let $\mathcal{T} \subset \mathcal{B}_{\Sigma_A}$ be a countably generated sigma-algebra, and let $\mathcal{D}$ be a countable generator of $\mathcal{T}$. We write $\{[x]_{\mathcal{T}} : x \in \Sigma_A\}$ for the space of atoms of $\mathcal{T}$, where $[x]_{\mathcal{T}} = \bigcap_{x \in D \in \mathcal{D}} D$ for every $x \in \Sigma_A$, and denote by $\{\mu_x^{\mathcal{T}} : x \in \Sigma_A\}$ a regular decomposition of $\mu$ over $\mathcal{T}$, i.e. a family of Borel probability measures on $\Sigma_A$ with the following properties:

(a) $\mu_x([x]_{\mathcal{T}}) = 1$ for every $x \in \Sigma_A$,
(b) the map $x \mapsto \mu_x^{\mathcal{T}}(B)$ from $\Sigma_A$ to $\mathbb{R}$ is $\mathcal{T}$-measurable for every Borel set $B \subset \Sigma_A$, and

$$\mu_x^{\mathcal{T}}(B) = E_\mu(1_B | \mathcal{T})(x) \quad \mu\text{-a.e.}$$

In view of (b) we may also assume without loss in generality that

(c) $\mu_{x'}^{\mathcal{T}} = \mu_x^{\mathcal{T}}$ for every $x \in \Sigma_A$ and $x' \in [x]_{\mathcal{T}}$.

Let $\mathcal{P}_0 = \{[0]_0, \ldots, [n-1]_0\}$ be the state (time-0) partition of $\Sigma_A$, and let $\mathcal{A} = \bigvee_{k \geq 0} \sigma^{-k}(\mathcal{P}_0) \subset \mathcal{B}_{\Sigma_A}$ be the sigma-algebra generated by all cylinder sets of the form $C = [i_0, \ldots, i_m]_0$, where $m \geq 0$ and $(i_0, \ldots, i_m) \in \{0, \ldots, n-1\}^{m+1}$. Then $\mathcal{A}^- = \sigma^{-1}(\mathcal{A}) \subset \mathcal{A}$, and the atoms of $\mathcal{A}$ and $\mathcal{A}^-$ are of the form

$$[x]_{\mathcal{A}} = \{x' \in \Sigma_A : x'_k = x_k \text{ for } k \geq 0\},$$

$$[x]_{\mathcal{A}^-} = \{x' \in \Sigma_A : x'_k = x_k \text{ for } k \geq 1\}$$

for every $x \in \Sigma_A$. The conditional information function $I_\mu(\mathcal{A}|\mathcal{A}^-)$ is given by

$$I_\mu(\mathcal{A}|\mathcal{A}^-)(x) = -\log E_\mu(1_{[x_0]}|\mathcal{A}^-)(x) = -\log \mu_x^{\mathcal{A}^-}([x]_{\mathcal{A}})$$

$$= \log \frac{d\mu_x^{\mathcal{A}}}{d\mu_x^{\sigma^{-1}(\mathcal{A})}}(x) \tag{6.2}$$

for every $x \in \Sigma_A$ (cf. (3.11) in [7]). The map

$$J = e^{-I_\mu(\mathcal{A}|\mathcal{A}^-)} \colon \Sigma_A \longmapsto [0,1] \subset \mathbb{R} \tag{6.3}$$

is obviously $\mathcal{A}$-measurable, and we claim that it is a function of the coordinates $x_0, x_1$ of every point $x \in \Sigma_A$.

In order to prove this claim we assume that $s = s_0 \ldots s_m$ and $s' = s'_0 \ldots s'_m$ are allowed strings with $s_0 = s'_0$, $s_1 = s'_1$ and $s_m = s'_m$, and such that the entries in $s$ and $s'$ are permutations of each other, and consider the map $V \in [S_A]$ which interchanges $s$ and $s'$, i.e.

$$(Vx)_k = \begin{cases} s'_k & \text{if } x \in [s]_0 \text{ and } 0 \leq k \leq m, \\ s_k & \text{if } x \in [s']_0 \text{ and } 0 \leq k \leq m, \\ x_k & \text{if } x \notin ([s]_0 \cup [s']_0) \text{ or } k > m \end{cases}$$

for every $x \in \Sigma_A$. Then $V(\mathcal{A}) = \mathcal{A}$, $V(\mathcal{A}^-) = \mathcal{A}^-$, and (6.2)–(6.3) show that

$$J(x) = J(Vx)$$

for every $x \in \Sigma_A$. The $\mathcal{A}$-measurability of $J$ allows us to view it as a map $J \colon \Sigma_A^+ \longmapsto \mathbb{R}$, and the argument in the preceding paragraph shows that this map is invariant under the subrelation of $S_A^+$ consisting of all $(x, x') \in S_A^+$ with $x_0 = x'_0$ and $x_1 = x'_1$. According to Lemma 6.1 this implies that $J$ is constant $\mu^+$-a.e. on each cylinder set $[i_0, i_1]_0 \subset \Sigma_A^+$, where $\mu^+$ is the projection of $\mu$ onto $\Sigma_A^+$, and thus a function of the coordinates $x_0, x_1$. This also proves our claim for the original map $J \colon \Sigma_A \longmapsto \mathbb{R}$.

Define a 0-1-matrix $A' = (A'(i, j), 0 \leq i, j \leq n-1)$ by setting (with the understanding that $0^0 = 0$)

$$A'(i, j) = \begin{cases} J(x_0, x_1)^0 \in \{0, 1\} & \text{if } x \in [i, j]_0 \subset \Sigma_A^+, \\ 0 & \text{if } [i, j]_0 = \varnothing \end{cases}.$$

Clearly $A'(i, j) \leq A(i, j)$ for every $(i, j) \in \{0, \ldots, n-1\}$. We claim that the SFT $\Sigma_{A'} \subset \Sigma_A$ defined by $A'$ has the following properties:

(d) $\Sigma_{A'}$ is equal to the (closed) support of $\mu$;

(e) if $R_{A'}^+$ is the Gibbs relation of the one-sided SFT $\Sigma_{A'}^+$, then

$$\frac{d\mu^+(x)}{d\mu^+(y)} = \prod_{k \geq 0} \frac{J(\sigma^k(x))}{J(\sigma^k(y))}$$

for every $(x, y) \in R_{A'}^+$;

(f) if $R_{A'}$ is the Gibbs relation of $\Sigma_{A'}^+$, then

$$\frac{d\mu^+(x)}{d\mu^+(y)} = \prod_{k \geq 0} \frac{J(\sigma^k(x))}{J(\sigma^k(y))}$$

for every $(x, y) \in R_{A'}$;

(g) $S_A(\Sigma_{A'}) = \Sigma_{A'}$.

Indeed, (e) is an immediate consequence of the definition of $J$ in terms of conditional measures, (f) follows from (e) and the shift-invariance of $\mu$, (d) follows from the fact that $\mu$ is the Gibbs measure of $\log J$ on $\Sigma_{A'}$, and (g) is a consequence of this as well as the $S_A$-invariance of $\mu$. Since $\mu$ is shift-invariant and $S_A$-ergodic, Remark 4.5 shows that the matrix $A'$ is aperiodic.

Next we prove that the one-sided measure $\mu^+$ is equivalent to an $S_A^+$-invariant measure $\nu^+$.

Renumber the symbols $0, \ldots, n-1$ of $A$, if necessary, choose an integer $n' \in \{1, \ldots, n-1\}$ such that first $n'$ rows and columns of $A'$ are nonzero and the remaining ones zero, and view $S = \{0, \ldots, n'-1\}$ as the alphabet of $\Sigma_{A'} \subset \Sigma_A \cap S^{\mathbb{Z}} \subset \{0, \ldots, n-1\}^{\mathbb{Z}}$. Put $P(i, j) = J(x_0, x_1)$ whenever $x = (x_k) \in \Sigma_{A'}$ satisfies that $x_0 = i, x_1 = j$. As $\mu$ is invariant under $S_A$, the conditions (f) and (g) above show that

$$\frac{d\mu(x)}{d\mu(y)} = \prod_{k \in \mathbb{Z}} \frac{P(x_k, x_{k+1})}{P(y_k, y_{k+1})} = 1 \tag{6.4}$$

whenever $(x, y) \in S_{A'} = S_A \cap (\Sigma_{A'} \times \Sigma_{A'})$. We set $S_{A'}^+ = S_A^+ \cap (\Sigma_{A'} \times \Sigma_{A'})$ and define the $S_{A'}^+$-equivalence classes of symbols of $\Sigma_{A'}^+ \subset \{0, \ldots, n'-1\}^{\mathbb{N}}$ as in Remark 4.4. For each such equivalence class $C \subset \{0, \ldots, n'-1\}$ choose and fix a symbol $i_C \in C$. Given any $x \in \Sigma_{A'}^+$, find the equivalence class $C$ with $x_0 \in C$, choose $y \in S_{A'}^+(x)$ with $y_0 = i_C$ and set

$$b(x) = \frac{d\mu^+(x)}{d\mu^+(y)}.$$

Then (from (6.4)) $b(x)$ is well-defined, i.e. independent of the choice of $y$, and only depends on the coordinate $x_0$ of $x$. For any $(x, y) \in S_{A'}^+$,

$$\frac{d\mu^+(x)}{d\mu^+(y)} = \frac{b(x)}{b(y)}.$$

Let $\nu^+$ be the unique probability measure on $\Sigma_{A'}^+$ which is a constant multiple of $b^{-1}\mu^+$. Then $\nu^+$ is Markov,

$$\frac{d\nu^+(x)}{d\nu^+(y)} = \frac{b(y)}{b(x)} \prod_{k \geq 0} \frac{P(x_k, x_{k+1})}{P(y_k, y_{k+1})}$$

for all $(x, y) \in R_{A'}^+$, and $\nu^+$ is invariant under $S_{A'}^+$ (which is equivalent to saying that $\nu^+$, regarded as a probability measure on $\Sigma_A^+$, is invariant under $S_A^+$).

Finally we show that $\mu$ can also be written as the Gibbs measure of a function of a single coordinate on $\Sigma_{A'} \subset \{0, \ldots, n'-1\}^{\mathbb{Z}}$.

Regard $b$ as a function on the alphabet $\{0, \ldots, n'-1\}$ of $\Sigma_{A'}$ by setting $b(i) = b(x)$ whenever $x_0 = i$, $i = 0, \ldots, n'-1$, and put $\bar{P}(i, j) = b(i)^{-1} P(i, j) b(j)$ for every $(i, j) \in \{0, \ldots, n'-1\}^2$ with $A'(i, j) = 1$. We claim that $\bar{P}(i, j) = \bar{P}(i, j')$ whenever $A'(i, j) = A'(i, j') = 1$ (or, equivalently, whenever $\mu([ij]_0)\mu([ij']_0) > 0$). As we have seen above,

$$\frac{d\nu^+(x)}{d\nu^+(y)} = \prod_{k \geq 0} \frac{\bar{P}(x_k, x_{k+1})}{\bar{P}(y_k, y_{k+1})}$$

for every $(x, y) \in R_{A'}^+$. Given $A'$-allowable 2-blocks $ij$ and $ij'$, Lemma 6.1 implies that there exists a pair $(x, y) \in S_{A'}^+ \cap ([ij]_0 \times [ij']_0) = S_A^+ \cap (\Sigma_{A'}^+ \times \Sigma_{A'}^+) \cap ([ij]_0 \times [ij']_0)$. As $\nu^+$ is $S_A^+$-invariant,

$$\frac{d\nu^+(x)}{d\nu^+(y)} = \prod_{k \geq 0} \frac{\bar{P}(x_k, x_{k+1})}{\bar{P}(y_k, y_{k+1})} = 1.$$

However, the pair $(\sigma(x), \sigma(y))$ also lies in $S_{A'}^+ = S_A^+ \cap (\Sigma_{A'}^+ \times \Sigma_{A'}^+)$, since $x_0 = y_0$, so that

$$\frac{d\nu^+(\sigma(x))}{d\nu^+(\sigma(y))} = \prod_{k \geq 1} \frac{\bar{P}(x_k, x_{k+1})}{\bar{P}(y_k, y_{k+1})} = 1.$$

By comparing the last two equations we have established that $\bar{P}(i, j) = \bar{P}(i, j')$, which proves our claim.

Let $\phi \colon \{0, \ldots, n'-1\} \longmapsto \mathbb{R}$ be the map satisfying that

$$\phi(x_0) = \log \bar{P}(x_0, x_1) = -I_\mu(x_0, x_1) - \log b(x_0) + \log b(x_1)$$

for every $x \in \Sigma_{A'}$. From condition (f) above we know that

$$\rho_\mu(x,y) = \sum_{k \in \mathbb{Z}} (\phi(x_k) - \phi(y_k))$$

for every $(x,y) \in R_{A'}$. In other words, $\mu$ is the unique Gibbs measure of $\phi$ on $\Sigma_{A'}$.

This completes the proof of the first part of the theorem, and the converse assertion is obvious from (3.4), the definition of $S_A$ and Theorem 3.3. □

*Remark 6.3.* If $\mu$ is atomic, invariant and ergodic under $S_A$, then $\mu$ has to be concentrated on a fixed point.

*Example 6.4. Symmetric measures on the Golden Mean SFT.* If $A = \left( \begin{smallmatrix} 1 & 1 \\ 1 & 0 \end{smallmatrix} \right)$, then Theorem 2.11 describes all the $S_A^+$-invariant probability measures on $\Sigma_A$. In interpreting this statement in the two-sided case one has to be a little careful (cf. (6.1)): the two points $\dots 01010101010 \dots$ of period two are both fixed under $S_A$; therefore each of them carries an $S_A$-invariant probability measure which is not shift-invariant. However, every *nonatomic* $S_A$-invariant and ergodic probability measure on $\Sigma_A$ is also shift-invariant.

*Example 6.5. A nonatomic, symmetric but not shift-invariant measure on an SFT.* If

$$A = \left( \begin{smallmatrix} 1 & 1 & 0 \\ 1 & 1 & 1 \\ 0 & 1 & 1 \end{smallmatrix} \right),$$

then the point $\dots 000012222 \dots$ is fixed under $S_A$, but not under the shift, and carries an atomic, $S_A$-invariant probability measure which is not shift-invariant, and whose orbit under the shift is infinite. In order to obtain a nonatomic example we set

$$A = \left( \begin{smallmatrix} 1 & 1 & 1 & 0 & 0 \\ 1 & 1 & 1 & 0 & 0 \\ 1 & 1 & 1 & 1 & 1 \\ 0 & 0 & 1 & 1 & 1 \\ 0 & 0 & 1 & 1 & 1 \end{smallmatrix} \right)$$

and consider the unique $S_A$-invariant probability measure on the set $\{x = (x_k) \in \Sigma_A : x_0 = 2$ and $x_{-k} \in \{0,1\}, x_k \in \{3,4\}$ for every $k > 0\}$.

*Example 6.6. An irreducible, aperiodic, $S_A$-invariant SFT $\Sigma_{A'} \subset \Sigma_A$ which is not determined by its alphabet.* Let

$$A = \left( \begin{smallmatrix} 1 & 1 & 1 \\ 0 & 0 & 1 \\ 1 & 0 & 0 \end{smallmatrix} \right), \qquad A' = \left( \begin{smallmatrix} 1 & 1 & 0 \\ 0 & 0 & 1 \\ 1 & 0 & 0 \end{smallmatrix} \right).$$

Then $\Sigma_{A'}$ has the same alphabet as $\Sigma_A$ and is irreducible, aperiodic and $S_A$-invariant. This example shows that the irreducible, aperiodic and $S_A$-invariant SFT's $\Sigma_{A'} \subset \Sigma_A$ which arise in Theorem 6.2 are not necessarily of the form $\Sigma_A \cap S^{\mathbb{Z}}$ for some subset $S \subset \{0, \dots, n-1\}$.

## 7. QUESTIONS AND REMARKS

7.1. **Weak mixing.** As in Section 1 we let $V = (V_k, k \geq 0)$ be a sequence of nonempty, finite, totally ordered sets, put $X_V = \prod_{k \geq 0} V_k$, assume that $Y \subset X_V$ is a nonempty, closed subset satisfying (M), write $T_Y$ for the adic transformation of $Y$, and define the $m$-weights $w_m(y)$, $y \in Y$, by (2.13). In the case of a full shift with an adic-invariant measure, $T_Y$ has an eigenfunction with eigenvalue $\zeta$ only if $\zeta^{w_m(y)} \to 1$ $\mu$-a.e. How general is this statement?

Conditions for existence of eigenfunctions in the stationary case were developed in [34] and [51]. We conjecture that even in the measure-preserving Bernoulli case on the full 2-shift the adic transformation $T_Y$ is weakly mixing: $\zeta^{w_m(y)} \to 1$ $\mu$-a.e. implies that $\zeta = 1$. It is not difficult to see that $T_Y$ is topologically weakly mixing. Also, because of the self-similar structure of Pascal's triangle modulo a prime [30, 36], $T_Y$ cannot have an eigenvalue (other than 1) that is a root of unity. If $\zeta$ is not a root of unity, perhaps $\zeta^{w_m(y)}$ is even almost surely uniformly distributed in the unit circle, or at least dense. Conceivably, it is not convergent to 1 down *every* path inside Pascal's triangle. Consider the skew product transformation $T(z_1, z_2, z_3, \dots) = (\zeta z_1, z_1 z_2, z_2 z_3, \dots)$ on the infinite torus. For $\zeta$ not a root of unity, it is well known that this is uniquely ergodic, and so, as we look at any fixed coordinate in the successive members of the orbit of any point, we see a uniformly distributed sequence. If we flip a coin with probabilities $\alpha$ and $1 - \alpha$ of heads and tails, and each time the coin comes up heads we shift our view one coordinate to the right, will we still see, with probability 1, a uniformly distributed, dense, or at least not convergent to 1 sequence?

7.2. **Super-$K$.** Let us say that a finite-state process $(\mathcal{P}, T)$, where $T$ is a measure-preserving or nonsingular transformation on a probability space $(X, \mu)$ and $\mathcal{P}$ is a finite measurable partition of $X$, is (one- or two-sided) *super-$K$* if the associated (dependent, transient) random walk on $\mathcal{P} \times \{0, 1, 2, \dots\}^{|\mathcal{P}|}$, which keeps track not only of which symbol appears (i.e., which atom of $\mathcal{P}$ is entered) at each time $n$, but also *how many times* each symbol has appeared up to time $n$, has (one- or two-sided, respectively) trivial tail. We have shown above that certain processes with Gibbs measures, including all mixing Markov processes, are two-sided super-$K$. Are there other natural examples, for example processes with the right uniform rate of mixing? Does every $K$-system have a super-$K$ generator?

7.3. **Invariant Measures.** In Section 6 we identified the shift- and $S_A$-invariant measures on subshifts of finite type. More generally, if $\psi$ is a continuous map on $\Sigma_A$ taking values in an arbitrary countable, discrete group $\mathcal{G}$ with finite conjugacy classes, what can be said about the set of $S_A^\psi$-invariant measures on $\Sigma_A$? It is not difficult to check that, if $\psi$ takes values in a *finite* group, then the only $S_A^\psi$-invariant measure is the measure of maximal entropy on $\Sigma_A$ (i.e. the unique $R_A$-invariant measure). More generally, if $\mathcal{G}$ is a countable, discrete group with finite conjugacy classes and $\mathcal{G}'$ a finite group, and if $\psi\colon \Sigma_A \longmapsto \mathcal{G}$ and $\psi'\colon \Sigma_A \longmapsto \mathcal{G}'$ are continuous maps, we denote by $(\psi, \psi')\colon \Sigma_A \longmapsto \mathcal{G} \times \mathcal{G}'$ the product map and leave it as an exercise to show that the set of nonatomic, shift-invariant measures which are invariant under $S_A^\psi$ coincides with the corresponding set for $S_A^{(\psi,\psi')}$.

7.4. **Asymptotic formulae.** Let $A = (A(i,j), 0 \le i, j \le n-1)$ be an irreducible and aperiodic 0-1-matrix, $P = (P(i,j), 0 \le i, j \le n-1)$ a stochastic matrix which is compatible with $A$, and $\bar{p}$ the probability vector satisfying $\bar{p}P = \bar{p}$. By (4.1), the shift-invariant Markov measure $\bar{\mu}_P$ in (4.3) is quasi-invariant under the relation $S_A$ appearing in Example 2.2, and is equivalent to the Gibbs measure $\mu_{\log P}$ (cf. (4.2)). We define $\psi_0\colon \Sigma_A^+ \longmapsto \mathbb{Z}^{\mathbb{Z}/n\mathbb{Z}}$ as in Examples 5.3 and 5.4 and set

$$\psi_0(0,x) = 0 \in \mathbb{Z}^{\mathbb{Z}/n\mathbb{Z}},$$

$$\psi_0(m,x) = \sum_{k=0}^{m-1} \psi_0(\sigma^k(x)),$$

$$[[s]]_m = \{x \in \Sigma_A^+ : \psi_0(m+1,x) = s\}$$

for every $m \ge 1$, $s \in \mathbb{Z}^{\mathbb{Z}/n\mathbb{Z}}$ and $x \in \Sigma_A$.

The following proposition, which gives necessary and sufficient conditions for a (quasi-invariant) Markov measure to be ergodic for the symmetric equivalence relation $S_A^+$ (equivalently the adic transformation $T^*_{\Sigma_A^+}$) on the SFT $\Sigma_A^+$, should be compared with Theorem 2.5.

**Proposition 7.1.** *The following conditions are equivalent.*

(1) *The measure $\bar{\mu}_P$ is ergodic under the adic transformation $T^*_{\Sigma_A^+}$ (or, equivalently, under the relation $S_A^+$);*

(2) *For every Borel set $B \subset \Sigma_A^+$,*

$$\lim_{m\to\infty} \frac{\bar{\mu}_P(B \cap [[\psi_0(m+1,x)]]_m \cap [x_m]_m)}{\bar{\mu}_P([[\psi_0(m+1,x)]]_m \cap [x_m]_m)}$$

$$\lim_{m\to\infty} \frac{\bar{\mu}_P(B \cap [[\psi_0(m+1,x)]]_m)}{\bar{\mu}_P([[\psi_0(m+1,x)]]_m)} = \bar{\mu}_P(B) \quad \bar{\mu}_P\text{-a.e.;}$$

(3) *For every* $i \in \{0, \dots n-1\}$,

$$\lim_{m \to \infty} \frac{\bar{\mu}_P([i]_0 \cap [[\psi_0(m+1,x)]]_m \cap [x_m]_m)}{\bar{\mu}_P([[\psi_0(m,x)]]_m \cap [x_m]_m)}$$

$$\lim_{m \to \infty} \frac{\bar{\mu}_P([i]_0 \cap [[\psi_0(m+1,x)]]_m)}{\bar{\mu}_P([[\psi_0(m,x)]]_m)} = \bar{p}(i) \ \ \bar{\mu}_P\text{-a.e.}$$

*Proof.* Put $R_{\Sigma_A^+}^{(*,m)} = \{(x,x') \in S_A^+ : x_k = x'_k \text{ for } k \geq m\}$, $m \geq 0$, write $\mathcal{B}_{\Sigma_A^+}^{(*,m)}$ for the sigma-algebra of $R_{\Sigma_A^+}^{(*,m)}$-saturated Borel subsets of $\Sigma_A^+$, and set

$$\mathcal{B}_{\Sigma_A^+}^* = \bigcap_{m \geq 0} \mathcal{B}_{\Sigma_A^+}^{(*,m)}.$$

For every $m \geq 0$, $R_A^{+(*,m)}$ is the smallest sigma-algebra with respect to which the maps $\psi_0(k, \cdot)$, $k \geq m$, are measurable.

As we saw in the proof of Theorem 4.3, there exists a finite partition $\mathcal{O}_1, \dots, \mathcal{O}_p$ of $\Sigma_A^+$ into closed and open subsets such that the restriction of $\bar{\mu}_P$ to each $\mathcal{O}_i$ is ergodic under $R_{\Sigma_A^+}^* \cap (\mathcal{O}_i \times \mathcal{O}_i)$. From the choice of the map $\bar{h}$ in (3.14) and Remark 4.4 it is furthermore clear that each of the sets $\mathcal{O}_i$ can be chosen as

$$\mathcal{O}_i = \bigcup_{j \in F_i} [j]_0,$$

where $F_1, \dots, F_p$ is the partition of the alphabet $\mathbb{Z}/n\mathbb{Z}$ into $R_{\Sigma_A^+}^*$-equivalence classes of symbols, and that

$$\mathcal{B}_{\Sigma_A^+}^* = \{\mathcal{O}_1, \dots, \mathcal{O}_p\} \pmod{\bar{\mu}_P}.$$

If $B \subset \Sigma_A^+$ is a Borel set, then the martingale convergence theorem implies that

$$\lim_{m \to \infty} E_{\bar{\mu}_P}(1_B | \mathcal{B}_{\Sigma_A^+}^{(*,m)}) = E_{\bar{\mu}_P}(1_B | \mathcal{B}_{\Sigma_A^+}^*)$$

$\bar{\mu}_P$-a.e. and in $L^1(\Sigma_A^+, \mathcal{B}_{\Sigma_A^+}, \bar{\mu}_P)$. In particular, if $\bar{\mu}_P$ is ergodic, then

$$\lim_{m \to \infty} \frac{\bar{\mu}_P(B \cap [[\psi_0(m+1,x)]]_m \cap [x_m]_m)}{\bar{\mu}_P([[\psi_0(m+1,x)]]_m \cap [x_m]_m)}$$

$$= \lim_{m \to \inf ty} \frac{1}{\bar{\mu}_P([[\psi_0(m+1,x)]]_m \cap [x_m]_m)}$$

$$\cdot \int E_{\bar{\mu}_P}(1_B | \mathcal{B}_{\Sigma_A^+}^{(*,m)}) \cdot 1_{[[\psi_0(m+1,x)]]_m \cap [x_m]_m} \, d\bar{\mu}_P$$

$$= \bar{\mu}_P(B) \ \ \bar{\mu}_P\text{-a.e.},$$

since $[[s]]_m \cap [i]_m \in \mathcal{B}^{(*,m)}_{\Sigma^+_A}$ for every $s \in \mathbb{Z}^{\{0,\ldots,sn-1\}}$ and $i \in \mathbb{Z}/n\mathbb{Z}$. The omission of the terms $[x_m]_m$ is equivalent to replacing $[[\psi_0(m+1,x)]]_m \cap [x_m]_m = [[\psi_0(m)]]_{m-1} \cap [[\psi_0(m+1,x)]]_m$ by $[[\psi_0(m+1,x)]]_m$ and does not affect the calculation. This shows that (1) implies (2) and hence (3).

Conversely, if $\bar{\mu}_P$ is nonergodic, the description of the sigma-algebra $\mathcal{B}^*_{\Sigma^+_A}$ shows that, for every $i \in \mathbb{Z}/n\mathbb{Z}$, $E_{\bar{\mu}_P}(1_{[i]_0}|\mathcal{B}^*_{\Sigma^+_A}) = 0$ on $\Sigma^+_A \smallsetminus \mathcal{O}_j$, where

$$\mathcal{O}_j = \bigcup\{[i']_0 \, i' \in \{0,\ldots,n-1\} \text{is } S^+_A\text{-equivalent to } i\}.$$

It follows that there exists, for every sufficiently large $m$, a set $C \in \mathcal{B}^{(*,m)}_{\Sigma^+_A}$ of the form $C = [[s_m]]_m \cap \cdots \cap [[s_{m+j}]]_{m+j}$ or, equivalently, of the form $C = [i]_m \cap [[s_{m+1}]]_{m+1} \cap \cdots \cap [[s_{m+j}]]_{m+j}$, with

$$\frac{\bar{\mu}_P(1_{[i]_0} \cap C)}{\bar{\mu}_P(C)} = \frac{1}{\bar{\mu}_P(C)} \int 1_C \cdot E_{\bar{\mu}_P}(1_{[i]_0}|\mathcal{B}^{(*,m)}_{\Sigma^+_A}) \, d\bar{\mu}_P$$
$$< \frac{\bar{\mu}_P([i]_0)}{2} = \frac{\bar{p}(i)}{2}.$$

Since the measure $\bar{\mu}_P$ is Markov,

$$\frac{\bar{\mu}_P([i]_0 \cap C)}{\bar{\mu}_P(C)} = \frac{\bar{\mu}_P([i]_0 \cap [[s_m]]_m \cap [i]_m)}{\bar{\mu}_P([[s_m]]_m \cap [i]_m)}.$$

From the two different ways of writing $C$ it is clear that the terms $[x_m]_m$ can be omitted in this calculation. We have proved that (3)—and hence (2)—cannot hold if $\bar{\mu}_P$ is nonergodic. □

Proposition 7.1 has combinatorial consequences whose direct proof appears quite difficult in all but the most elementary cases. Let $A = \left(\begin{smallmatrix} 1 & 1 \\ 1 & 1 \end{smallmatrix}\right)$, $\Sigma^+_A = \Sigma^+_2$, and let $P = \left(\begin{smallmatrix} a & b \\ c & d \end{smallmatrix}\right)$ be a stochastic matrix which is compatible with $A$, i.e. which satisfies that $abcd > 0$. We denote by $\mu = \bar{\mu}_P$ the shift-invariant probability measure defined in (3.1), where $\bar{p}$ is the probability vector satisfying $\bar{p}P = \bar{p}$. Since $\mu = \bar{\mu}_P$ is equivalent to the Gibbs measure $\mu_{\log P}$, Theorem 4.3 and Remark 4.4 implies that $\mu$ is ergodic under $S^+_2$, so that $\mu$ must satisfy the equivalent conditions (2) and (3) in Proposition 7.1.

For $i, j \in \{0,1\}$, $m \geq 1$ and $s = (s_0, s_1)$ with $s_0 + s_1 = m+1$, let

$$Q_{i,j}(s) = \frac{\mu([i]_0 \cap [[(s)]]_m \cap [j]_m)}{\bar{p}(i)}.$$

If $C = [i_0, \ldots, i_{m'}]_0 \subset \Sigma_2^+$ with $\sum_{j=0}^{m'} i_j = s_1'$, $s_0' = m' + 1 - s_1'$, and $s' = (s_0', s_1')$, then

$$\frac{\mu(C \cap [[s]]_m \cap [j]_m)}{\mu(C)} = Q_{i_{m'}j}(s - s'),$$

and Proposition 7.1 implies that

$$\lim_{m \to \infty} \frac{Q_{i_{m'}x_m}(\psi_0(m+1, x) - s')}{Q_{j_{m''}x_m}(\psi_0(m+1, x) - s'')} = 1 \tag{7.1}$$

for $\mu$-a.e. $x \in \Sigma_2^+$, and for any $s', s'' \in \mathbb{N}^2$. In other words, Proposition 7.1 is a statement about 'amnesia' concerning the initial coordinates of the points $x \in \Sigma_2^+$. In this special case the quantities $Q_{ij}(s)$ can be determined explicitly (as sums of products of binomial coefficients and other factors), but the verification of the convergence

$$\lim_{m \to \infty} \frac{Q_{ij}(s^{(m)} - s')}{Q_{i'j}(s^{(m)} - s'')} = 1$$

in (7.1) for suitable sequences $(s^{(m)} = (s_0^{(m)}, s_1^{(m)}), m \geq 1)$ with $s_i^{(m)} \to \infty$ for $i = 0, 1$ looks forbidding.

7.5. **Interval splitting.** There should be some consequences for interval splitting of the ergodicity of Markov measures for the Pascal adic on the full shift. Perhaps if intervals are split in different proportions depending on whether they arose as the left or right half of a previously-split interval the resulting division points would still be uniformly distributed. Weak mixing of the Pascal adic might imply uniform distribution in a rectangle if two intervals are split simultaneously.

7.6. **Dynamics.** What are the joinings of the adic transformation on the full 2-shift with two different Bernoulli measures? Are these adic systems disjoint, or at least not isomorphic? A very general problem is to describe the dynamics of the measures found in our above theorem, especially in case the equivalence relation $S_A^{\psi+}$ consists of the orbit relation for a single nonsingular transformation.

7.7. **Possible applications.** Subshifts of finite type are important for the design of actual communication systems, especially in magnetic recording. In biological or materials science applications (say polymer building) SFT's might arise from something like the momentary disabling of a receptor: receipt of a certain symbol could make the receptor momentarily sensitive only to certain other symbols. If we want to model a signal recorded with constraints of the kind imposed by membership in a subshift of finite type (for example to record efficiently on a disk that has already been used), in the absence of further information it might be reasonable to assume that the statistics of the signal are given by a measure with some of the symmetries discussed above. For example, as in de Finetti's motivation for exchangeability (the distribution of repeated samples should be independent of the order in which they are drawn), perhaps it is natural to assume that cylinder sets in the SFT which map to one another under permutations of finitely many coordinates (or are symmetric in some other respect) have equal probabilities. For applications, the case of higher-dimensional actions needs serious development. The general situations connected with $S_A^{\psi+}$-invariance (which can keep track of accumulated symbol or transition counts or derangement effects) might reflect a cost associated with sending or receiving certain signals, or a hysteresis or memory in materials. In image enhancement and pattern generation and recognition, symmetric Gibbs measures could account for the presence of texture, like bands or dapples [1, 13]. Knowing their dynamics, for example spectral properties, could help to filter out such background by using appropriate transforms, to code signals better by taking the structure into account, or to detect boundaries between regions with different textures. The presence of discrete spectrum might reflect subtle rhythms (periodic or aperiodic regularities different from spatial regularities, which are connected with the dynamics of the associated shift transformations). Weak mixing would be tantamount to the lack of any such almost periodic structure and might indicate resistance to filtration by Fourier methods.

DEPARTMENT OF MATHEMATICS, CB 3250, PHILLIPS HALL, UNIVERSITY OF NORTH CAROLINA, CHAPEL HILL, NC 27599 USA
*E-mail address*: `petersen@math.unc.edu`

DEPARTMENT OF MATHEMATICS, UNIVERSITY OF VIENNA, VIENNA, AUSTRIA
*E-mail address*: `klaus.schmidt@univie.ac.at`